\documentclass[onecolumn]{IEEEtran}
\usepackage{cite}
\usepackage{amsmath,amssymb,amsfonts}
\usepackage{algorithmic}
\usepackage{graphicx}
\usepackage{epstopdf}
\usepackage{textcomp}
\usepackage{multirow}
\usepackage{booktabs}
\usepackage{subfigure}
\usepackage{stfloats}
\usepackage{color}
\usepackage{bm}
\usepackage{array}
\usepackage{float}
\usepackage{color}
\usepackage{threeparttable}
\usepackage{booktabs}
\usepackage{tablefootnote}
\def\BibTeX{{\rm B\kern-.05em{\sc i\kern-.025em b}\kern-.08em
		T\kern-.1667em\lower.7ex\hbox{E}\kern-.125emX}}

\begin{document}
				
	\title{A Hybrid SIE-PDE Formulation Without Boundary Condition Requirement for Transverse Magnetic Electromagnetic Analysis}
	\author{\IEEEauthorblockN{Aipeng Sun, \IEEEmembership{Graduate Student Member, IEEE}, Zekun Zhu, \IEEEmembership{Graduate Student Member, IEEE}, \\ Shunchuan Yang, \IEEEmembership{Member, IEEE}, and Zhizhang (David) Chen, \IEEEmembership{Fellow, IEEE} \\}
		
		\thanks{Manuscript received xxx; revised xxx.}
		\thanks{This work was supported in part by the National Natural Science Foundation
			of China through Grant 61801010, 62071125, 61671257, in part by Pre-Research Project
			through Grant J2019-VIII-0009-0170 and Fundamental Research Funds for
			the Central Universities. {\it {(Corresponding author: Shunchuan Yang.)}}
			
			 A. Sun and Z. Zhu are with the School of Electronic and Information Engineering, Beihang University, Beijing, 100083, China (e-mail: sunaipeng1997@buaa.edu.cn, zekunzhu@buaa.edu.cn).
			
			S. Yang is with the Research Institute for Frontier Science and the School of Electronic and Information Engineering, Beihang University, Beijing, 100083, China (e-mail: scyang@buaa.edu.cn).
			
			Z. Chen is currently with the College of Physics and Information Engineering, Fuzhou University, Fuzhou, Fujian. P. R. China, on leave from the Department of Electrical and Computer Engineering, Dalhousie University, Halifax, Nova Scotia, Canada B3H 4R2 (email: zz.chen@ieee.org).
		}
	}
	
	\maketitle
	
	\begin{abstract}	
	A hybrid surface integral equation partial differential equation (SIE-PDE) formulation without the boundary condition requirement is proposed to solve the transverse magnetic (TM) electromagnetic problems. In the proposed formulation, the computational domain is decomposed into two \emph{overlapping} domains: the SIE and PDE domains. In the SIE domain, complex structures with piecewise homogeneous media, e.g., highly conductive media, are included. An equivalent model for those structures is constructed by replacing them with the background medium and introducing a surface equivalent electric current density on an enclosed boundary to represent their electromagnetic effects. The remaining computational domain and homogeneous background medium replaced domain consist of the PDE domain, in which inhomogeneous or non-isotropic media are included. Through combining the surface equivalent electric current density and the inhomogeneous Helmholtz equation, a hybrid SIE-PDE formulation is derived. It requires no boundary conditions, and is mathematically equivalent to the original physical model. Through careful construction of basis functions to expand electric fields and the equivalent current density, the discretized formulation is made compatible with the SIE and PDE domain interface. The accuracy and efficiency are validated through two numerical examples. Results show that the proposed SIE-PDE formulation can obtain accurate results, and significant performance improvements in terms of CPU time and memory consumption compared with the FEM are achieved.
	\end{abstract}
	
	\begin{IEEEkeywords}
		Finite element method, hybrid methods, method of moment, partial differential equation, surface equivalence theorem, surface integral equation
	\end{IEEEkeywords}
	
	\section{Introduction}
	Hybrid formulations show significant potential in solving challenging problems in electrical engineering fields, like scattering from multiscale and electrically large objects [\citen{PengNon-DDM2013}], parameters extraction in integrated circuits [\citen{GedneyFEM/MoM1992}] since they can inherit the merits and overcome possible disadvantages of each formulation. They have attracted much attention in the last two decades, and a number of formulations have been developed. Those formulations can be categorized into time-domain methods, like the hybrid finite-difference time-domain finite-element time-domain (FDTD-FETD) method [\citen{TeiFDTD-FETD2007}][\citen{YeungFDTD-FETD1999}], the hybrid spectral element time-domain finite-element time-domain (SETD/FETD) method [\citen{SunFDTD-SETD2019}], the unconditionally stable finite-element time-domain (US-FETD) method [\citen{FanUS-FETD2020}], and the finite-element time-domain generalized scattering matrix (FETD-GSM) method [\citen{ZhangFETD-GSM2017}], and frequency-domain methods, like the method of moment mode-matching (MoM/MM) method [\citen{AydoMoM/MM2019}][\citen{CabaMoM/MM2011}], the method of moment physical-optics (MoM-PO) method [\citen{ChenMoM-PO2007}-\citen{GongMoM-PO2006}], and the finite-element
	method/method of moment (FEM/MoM) [\citen{IlicFEM-MoM_ante_2009}-\citen{JiFEM/MoM2002}], which can efficiently solve various electromagnetic problems. In this paper, we mainly focus on those in the frequency domain, especially the hybrid FEM-MoM formulation.
	
	One main group of hybrid formulations in the frequency domain is to combine different types of integral equations, like the MoM/PO method [\citen{ChenMoM-PO2007}-\citen{GongMoM-PO2006}], the method of moment uniform geometrical theory of diffraction (MoM/UTD) method [\citen{TapMoM-UTD2005}][\citen{LiuMoM-UTD2010}], to solve the scattering problems induced by the electrically large and multiscale objects, and extract the electrical parameters from integrated circuits. However, the formulations described are incapable of handling objects with inhomogeneous or non-isotropic media since the surface equivalence theorem is not applicable. To mitigate this issue, another group of hybrid formulations combines the finite element and boundary integral formulations, like the FEM/MoM [\citen{IlicFEM-MoM_ante_2009}-\citen{JiFEM/MoM2002}], the finite element boundary integral 
	(FE-BI) formulation [\citen{EibertFE/BI1999}][\citen{ShengFE/BI2002}]. In those hybrid FE-BI formulations, the surface integral equation (SIE) formulation is often used to mimic the absorbing boundary conditions and truncate the computational domain in the FEM. Typical applications are the scalar FE-BI [\citen{MeyerFE-BI1994}],and its symmetrical version [\citen{ZhaoFEM-BEM2006}][\citen{YangFE-BI-MLFMA2013}], which are proposed to improve the efficiency, accuracy, and convergence properties. 
	
	In [\citen{Guanmultisolver2017}][\citen{Guanmultisolver2016}], the SIE formulation is used to model the electrically large structures with homogenous media, and the FEM is for geometrically fine structures or inhomogeneous media. In those formulations, the computational domain is decomposed into two \emph{non-overlapping} domains: one for the boundary element method (BEM) formulation and the other for the FEM formulation. A boundary condition, e.g., the second kind of boundary condition [\citen{Jin2015FEM}, Ch, 1, pp. 42-46], or the transmission condition [\citen{Jin2015FEM}, Ch. 14, pp. 1066-1083], is required to couple the FEM and BEM formulations on the interface of the two domains. Certain assumptions are required to ensure that electromagnetic fields meet those boundary conditions. However, such assumptions are not always satisfied in the practical physical fields. Therefore, they may suffer from the accuracy issue for near fields calculation and are mainly used to calculate far-fields [\citen{Guanmultisolver2017}][\citen{Guanmultisolver2016}]. 
	
	In [\citen{DodigFEM-BEM2021}], a hybrid BEM/FEM formulation with the edge element is proposed to solve time harmonic electromagnetic scattering problems. With the same type of elements as well as the same edge basis functions, the two formulation are coupled without any additional boundary conditions. On this basis, accurate near and far fields can be obtained. However, this hybrid formulation is still constrained by the FEM formulation when processing with the interior region. For example, when dealing with highly conductive media or geometrically fine structures, extremely fine meshes are required which will greatly increase the number of unknowns and decrease the efficiency.
	
	In essence, these hybrid formulations can be cast into the domain decomposition methods (DDMs), which implies that different formulations can be applied in the corresponding domains according to their geometrical and material properties [\citen{LeeDDM2005}-\citen{ZhaoDDM2008}]. Especially, one embedded domain decomposition method (DDM) is proposed for signal integrity design and electromagnetic scattering analysis in [\citen{LueDDM2018}][\citen{LuSI-DDM2018}]. The computational domain is decomposed into two \emph{overlapping} domains: mesh coarse and fine domains. The geometrically fine structures are placed in the fine mesh domain and moved away from the computational domain. Therefore, the remaining computational domain does not include the geometrically fine structures, and coarse mesh can be used. The transmission condition is used to couple fields between the two domains. It shows significant performance improvements in terms of CPU time and memory consumption. However, it may need several iterations to solve the problems and suffer from efficiency issues if the boundary condition is not properly handled [\citen{LueDDM2018}][\citen{LuSI-DDM2018}]. 
	
	To mitigate the above accuracy and efficiency issues in the existing hybrid methods, a novel hybrid SIE-PDE formulation without boundary condition requirement is proposed for electromagnetic analysis. This greatly expands our previous work in [\citen{SunSIE-PDE2021}]. The computational domain is first decomposed into two \emph{overlapping} domains like the embedded DDM. One domain is for the SIE formulation, which is bounded by an enclosed contour, and the other domain is for the PDE formulation, which includes the whole computational domain. In the SIE domain, an equivalent model is derived based on the surface equivalent theorem [\citen{LOVE}, Ch. 12, pp. 653-658]. The original structures are replaced by the background medium, and the surface equivalent electric current density is introduced on the boundary to keep fields unchanged in the exterior region. In the PDE domain, the inhomogeneous Helmholtz formulation incorporated with the surface equivalent electric current density on the boundary of the SIE domain is used to model inhomogeneous or non-isotropic media. Therefore, the electromagnetic fields in the whole computational domain can be formulated through a hybrid SIE-PDE formulation. Complex structures and computationally challenging media, like the highly conductive media, are replaced by the background medium, which can be easily handled with coarse meshes. Therefore, the proposed SIE-PDE formulation can significantly improve the efficiency over the traditional PDE formulations. The proposed formulation is mathematically equivalent to the original physical model, and requires no additional boundary conditions. Therefore, accurate near and far-fields can be obtained. 
	
	This paper is organized as follows. In Section II, configurations and preliminary notations are defined, and procedures for the proposed hybrid SIE-PDE formulation are illustrated. In Section III, the hybrid SIE-PDE formulation is proposed. In Section IV, detailed implementations for the proposed hybrid SIE-PDE formulation are shown. In addition, fields calculation in the whole computational domain is also discussed. In Section V, two numerical examples are carried out to validate its accuracy and efficiency. At last, we draw some conclusions in Section VI.
	
	\section{The Configurations and Preliminaries}
\begin{figure}[htbp]
		\begin{minipage}[t]{0.5\textwidth}
			\centering
			\centerline{\includegraphics[scale=0.24]{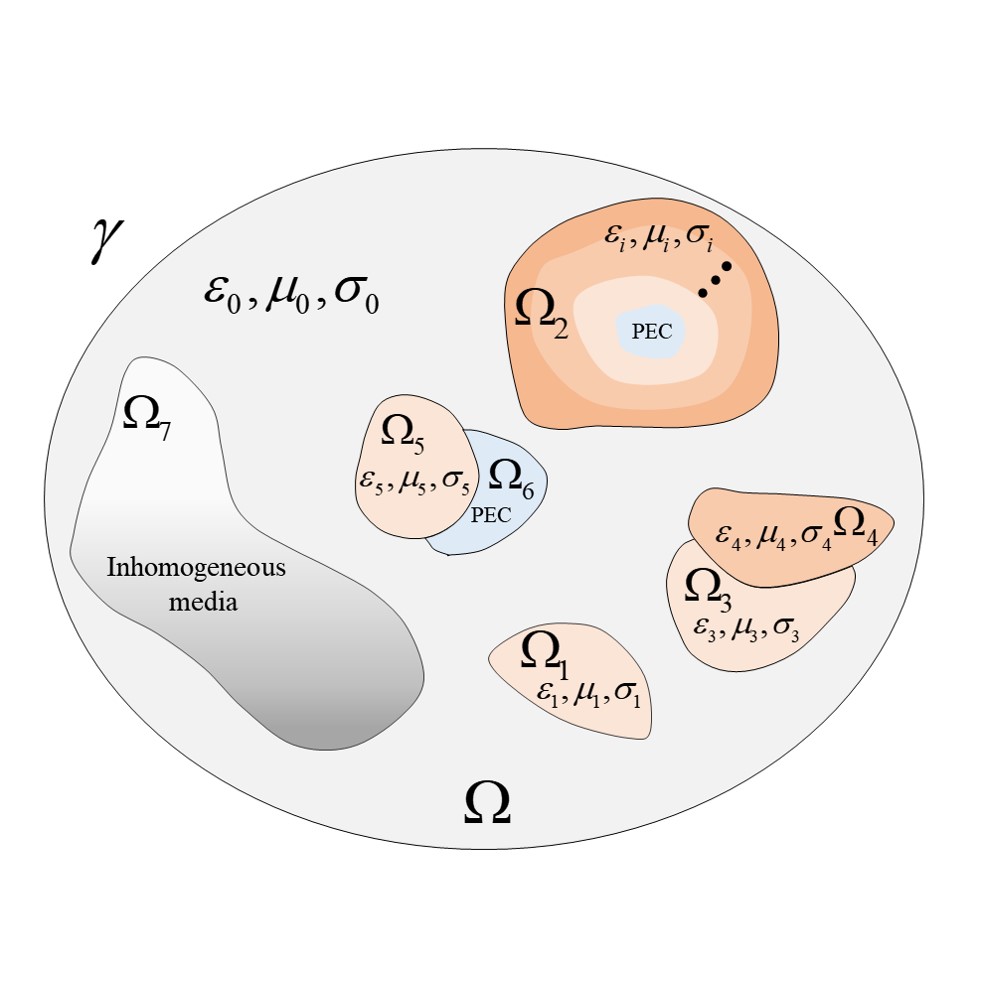}}
			\centering
			\centerline{(a)}
		\end{minipage}
		\begin{minipage}[t]{0.5\textwidth}
			\centering
			\centerline{\includegraphics[scale=0.24]{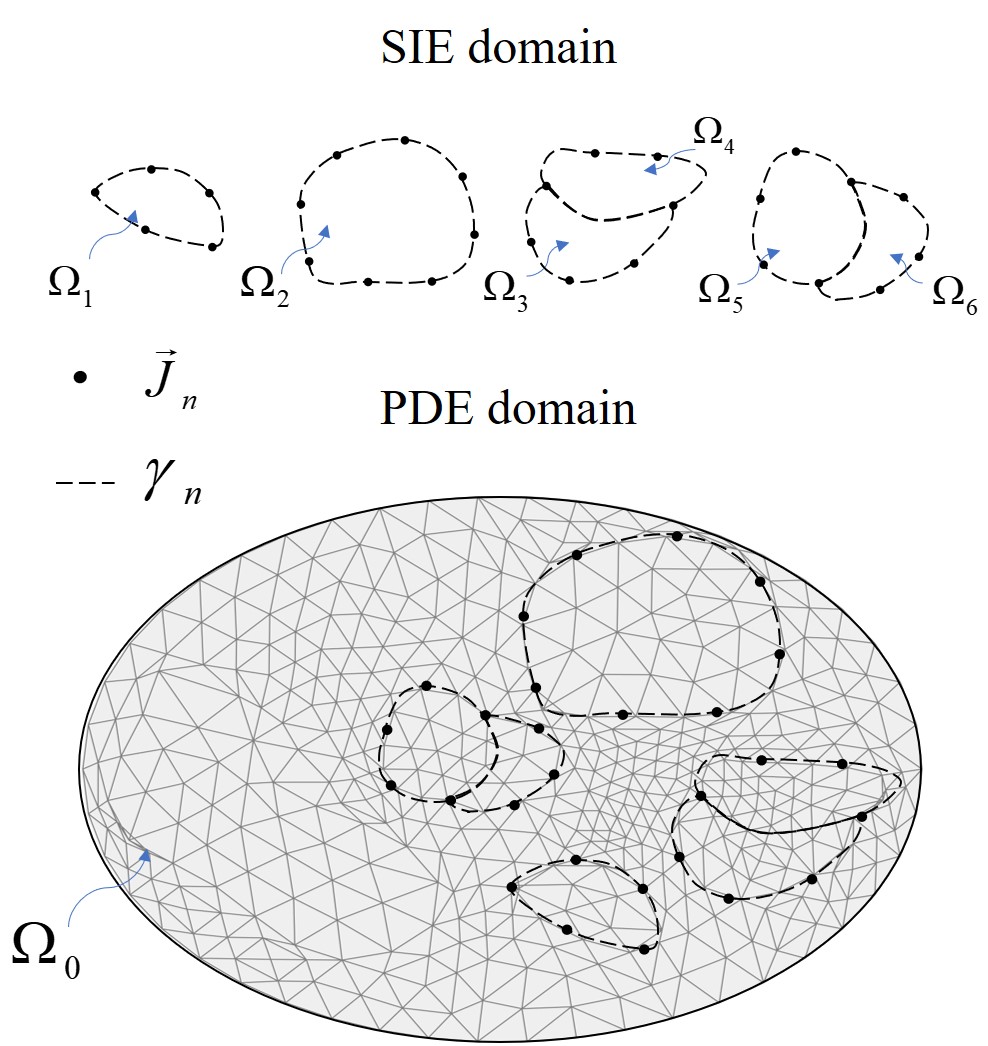}}
			\centering
			\centerline{(b)}
		\end{minipage}
		\caption{(a) The general scenario including complex structures with homogeneous and inhomogeneous media, and (b) the \emph{overlapping} SIE and PDE domains for the proposed hybrid SIE-PDE formulation. }
		\label{general_model}
	\end{figure}

	In this paper, as shown in Fig. \ref{general_model}(a), a general scenario with the two-dimensional (2D) TM scattering problem induced by objects with homogeneous and inhomogeneous media is considered to demonstrate the proposed SIE-PDE formulation. The constant parameters of those objects are the permittivity ${\varepsilon _i}$, the permeability ${\mu _i}$ and the conductivity ${\sigma _i}$, respectively, and its corresponding boundary ${\gamma _i}$. The background medium is air with constant parameters of ${\varepsilon _0}$, ${\mu _0}$, ${\sigma _0}$, and the whole computational domain is truncated by the outermost boundary $\gamma $. 
	
	A symbol with a subscript, e.g. ${E_i}$, denotes a continuous quantity on ${\gamma _i}$. A hollow symbol with a subscript, e.g. ${\mathbb{P}_i}$, denotes a matrix on $\gamma_i$, and a bold symbol with a subscript, e.g. ${\mathbf{E}_i}$, denotes a column vector on ${\gamma _i}$. Quantities with $\,\widehat{}\,\,$, e.g. $\widehat{E}$, are used for the equivalent model.
	
	The proposed hybrid SIE-PDE formulation needs the following three steps:
	\begin{enumerate}
		\item An appropriate boundary ${\gamma _n}$ is selected, and the computational domain is decomposed into two \emph{overlapping} domains: the SIE domain, which includes complex structures or computationally challenging media, and the PDE domain, in which inhomogeneous and non-isotropic media are filled as shown in Fig. \ref{general_model}. 
		
		\item In the SIE domain, an equivalent model is constructed through applying the surface equivalence theorem, and single surface electric current density $J_n$ is derived and enforced on ${\gamma _n}$.
		
		\item In the PDE domain, which includes the remaining computational domain and the background medium replaced domain, the inhomogeneous Helmholtz equation in terms of electric fields is used to model inhomogeneous and non-isotropic media. In addition, the surface equivalent electric current density on ${\gamma _n}$ is coupled into the Helmholtz equation to present the electromagnetic effects in the SIE domain. Then, the hybrid SIE-PDE formulation without the boundary condition requirement is derived.
	\end{enumerate}
	
	We reported the preliminary idea in [\citen{SunSIE-PDE2021}]. Detailed formulations and discussion will be presented in the following sections. We first derive the equivalent model incorporated with the single electric current density, and then derive the hybrid SIE-PDE formulation in terms of electric fields in the continuous physical domain.
	
	\section{The Proposed Hybrid SIE-PDE Formulation}
	\subsection{The Equivalent Model with the Single Electric Current Density for Objects with Piecewise Homogeneous Media in the SIE Domain}
\begin{figure}
		\begin{minipage}[t]{0.5\linewidth}
			\centerline{\includegraphics[scale=0.17]{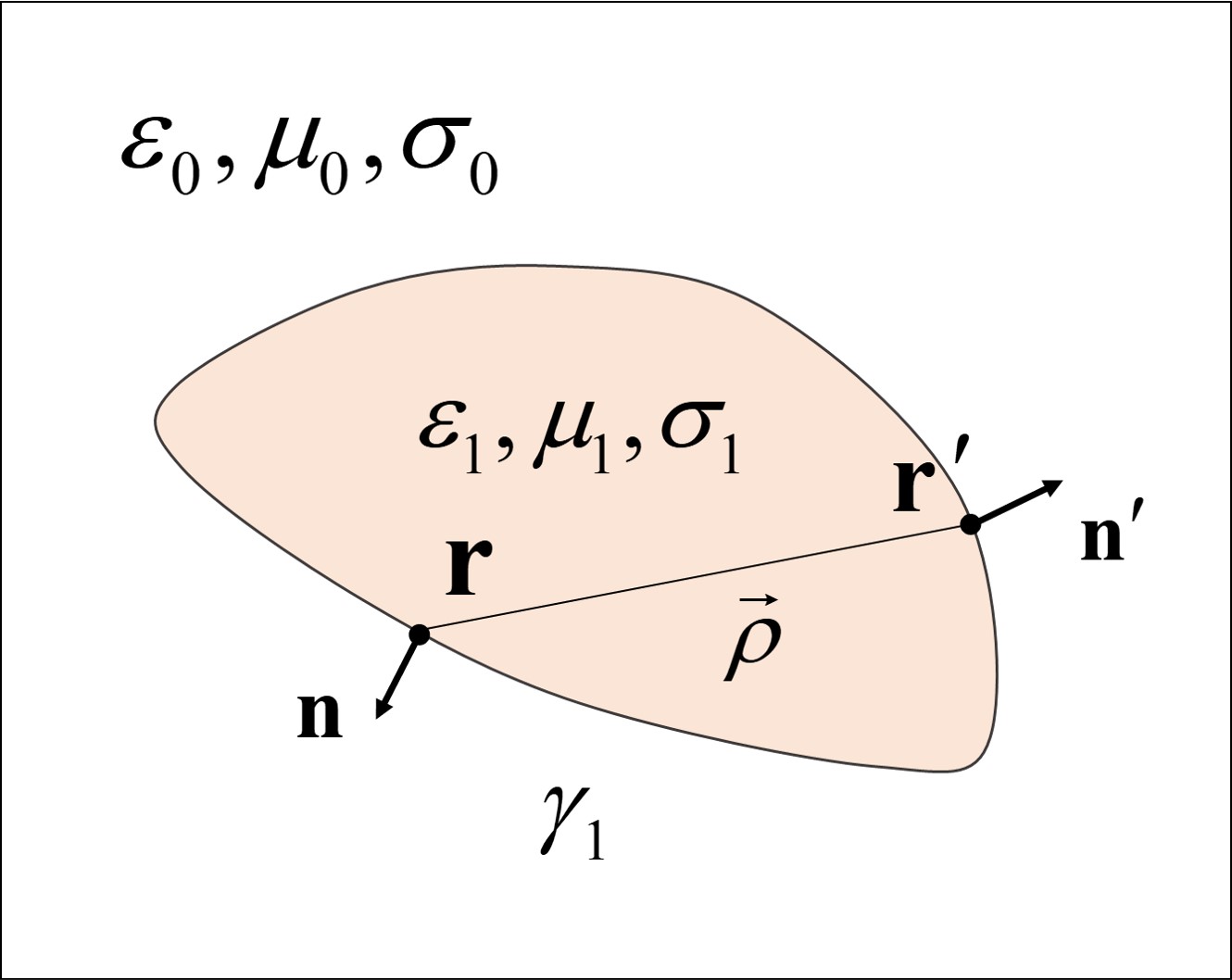}}
			\centerline{(a)}
		\end{minipage}
		\hfill
		\begin{minipage}[t]{0.5\linewidth}
			\centerline{\includegraphics[scale=0.17]{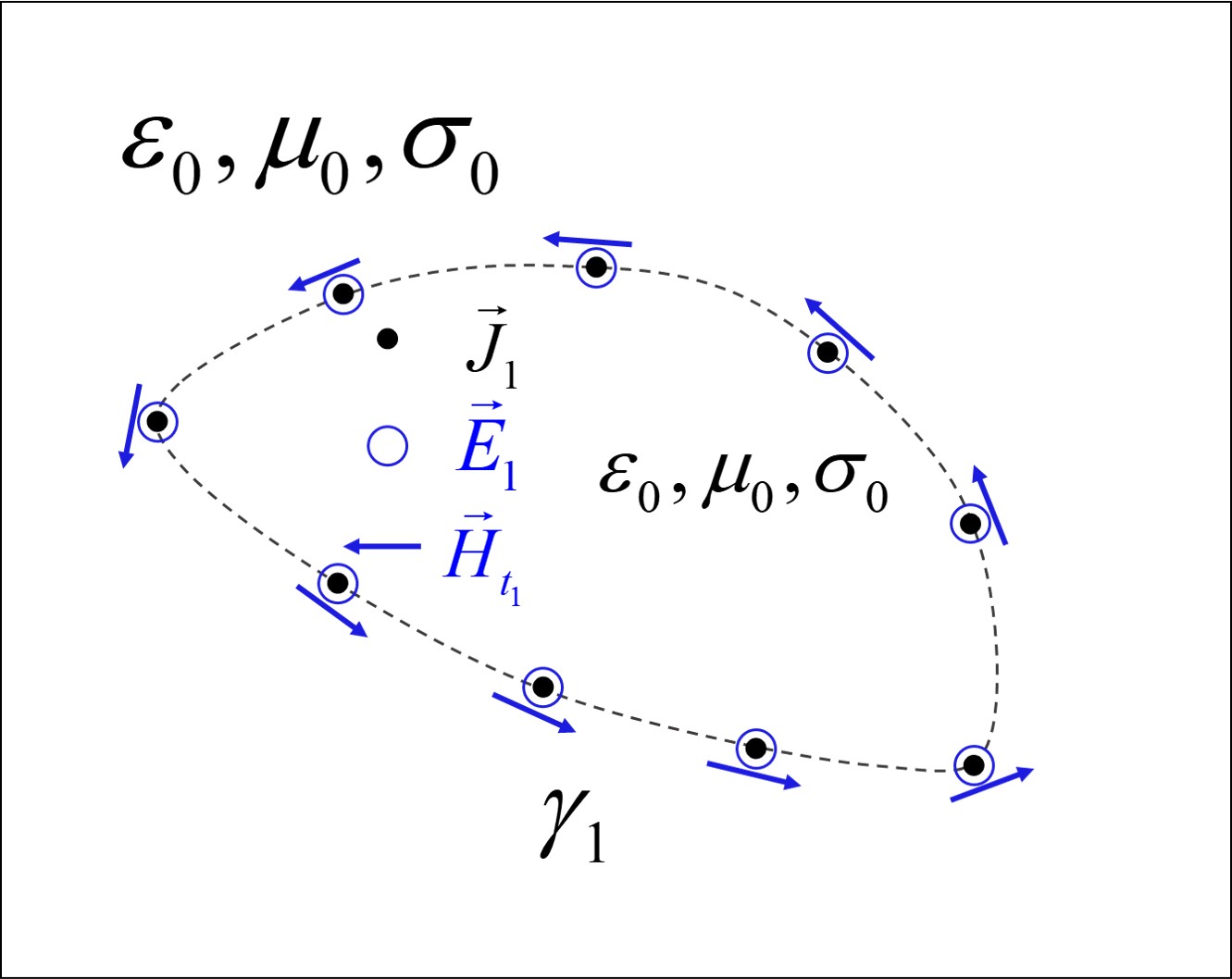}}
			\centerline{(b)}
		\end{minipage}
		\caption{(a) The original model for a single penetrable object, and (b) the equivalent model.}
		\label{single_model}
	\end{figure}

	We first illustrate the detailed derivation of the equivalent model with the single electric current density for objects with piecewise homogeneous media. To make the derivation concise, a single penetrable object with the permittivity ${\varepsilon _1}$, permeability ${\mu _1}$, and conductivity ${\sigma _1}$ is first considered, as shown in Fig. \ref{single_model}(a). Its boundary is denoted as ${\gamma _1}$. According to the surface equivalence theorem [\citen{LOVE}, Ch. 12, pp. 653-658], an equivalent model can be obtained, in which the object is replaced by its surrounding medium, and a surface equivalent electric current density is introduced on $\gamma_1$. To further enforce the electric fields in the original and equivalent configurations equal to each other on ${\gamma _1}$, a single-source (SS) formulation with only the electric current density on ${\gamma _1}$ is obtained in Fig. \ref{single_model}(b) [\citen{Zhou2021embedded}].
	
	Without loss of generality, the penetrable object does not include any sources. Therefore, the electric fields inside the penetrable object must satisfy the following homogeneous scalar Helmholtz equation
	
	\begin{equation} \label{hoHe}
		{\nabla ^2}{E_1}{\rm{ + }}{k^2}{E_1} = 0,
	\end{equation}
	subject to the boundary condition	
	\begin{equation} \label{BCE1}
		{\left. {{E_1}\left( {\mathbf{r}} \right)} \right|_{\mathbf{r} \in {\gamma _1}}} = {\left. {{{\widehat E}_1}\left( {\mathbf{r}} \right)} \right|_{\mathbf{r} \in {\gamma _1}}},
	\end{equation}
	where ${E_1}$, ${\widehat E_1}$ denote the electric fields inside ${\gamma _1}$ for the original and equivalent models, respectively, and ${\left. {{E_1}\left( {\mathbf{r}} \right)} \right|_{\mathbf{r} \in {\gamma _1}}}$, ${\left. {{{\widehat E}_1}\left( {\mathbf{r}} \right)} \right|_{\mathbf{r} \in {\gamma _1}}}$ denote their values on the inner side of ${\gamma _1}$.
	
	(\ref{hoHe}) can be solved through the contour integral method [\citen{Okoshi1985}]. Then, the electric fields inside ${\gamma _1}$ can be expressed in terms of ${E_1}$ and its normal derivative on ${\gamma _1}$ as	
	\begin{equation} \label{CIM}
		T{E_1}\!\left( {\mathbf{r}} \right) = \oint_{{\gamma _1}} {\left[ {{G_1}\left( 	{\mathbf{r},\mathbf{r}'} \right)\!\frac{{\partial {E_1}\left( {\mathbf{r}'} \right)}}{{\partial n'}} - {E_1}\left( {\mathbf{r}'} \right)\!\frac{{\partial {G_1}\left( {\mathbf{r},\mathbf{r}'} \right)}}{{\partial n'}}} \right]} d\mathbf{r}',
	\end{equation}		
	where $T = 1/2$, when the source point $\mathbf{r}'$ and the observation point $\mathbf{r}$  are located on the same boundary, otherwise, $T = 1$. ${G_1}\left({\mathbf{r},\mathbf{r}'}\right)$ is the Green’s function, which is given by ${G_1}\left( {\mathbf{r},\mathbf{r}'} \right) =-jH_0^{\left( 2 \right)}\left( {{k_1}\rho } \right)/4$, where $j = \sqrt {- 1}$, $\rho  = \left| {\mathbf{r}-\mathbf{r}'} \right|$, ${k_1}$ denotes the wave number in the penetrable object, and $H_0^{\left( 2 \right)}\left(  \cdot  \right)$ denotes the zeroth-order Hankel function of the second kind. $\mathbf{n}$, $\mathbf{n}'$ denote the unit normal vector pointing outward ${\gamma _1}$ at the observation point $\mathbf{r}$ and the source point $\mathbf{r}'$, respectively. In addition, the tangential magnetic field on ${\gamma _1}$, which is related to the electric field through the Poincare-Steklov operator [\citen{Patel2017SS-SIE}], can be written as
	
	\begin{equation} \label{TH1_o}
		{{H_{t_1}}\left( {\mathbf{r}} \right)}  = \frac{1}{{j\omega {\mu _1}}}{\left. {\frac{{\partial {E_1}\left( {\mathbf{r}} \right)}}{{\partial n}}} \right|_{\mathbf{r} \in {\gamma _1}}},
	\end{equation}
	where ${\mu _1}$ is the permeability of the penetrable object. We introduce the surface admittance operator (SAO) $\mathcal{Y}$ [\citen{Knockaert2008}], which relates the normal derivative of electric fields to electric fields on ${\gamma _1}$ . Therefore, it can be expressed as

	\begin{equation} \label{y1_define}
		{\left. {\frac{{\partial {E_1}\left( {\mathbf{r}} \right)}}{{\partial n}}} \right|_{\mathbf{r} \in {\gamma _1}}} = {\mathcal{Y}_1}{\left. {{E_1}} \right|_{\mathbf{r} \in {\gamma _1}}}.
	\end{equation}
	By substituting (\ref{y1_define}) into (\ref{CIM}), we have
	\begin{equation} \label{CIM_y1}
		T{E_1}\left( {\mathbf{r}} \right) = \oint_{{\gamma _1}} {\left[ {{G_1}\left( {\mathbf{r} ,\mathbf{r}'} \right){\mathcal{Y}_1} - \frac{{\partial {G_1}\left( {\mathbf{r},\mathbf{r}'} \right)}}{{\partial n'}}} \right]} {E_1}\left( {\mathbf{r}'} \right)d\mathbf{r}'. 
	\end{equation}
	To make the derivation more concise, we borrow some notations from [\citen{Knockaert2008}]. Then, (\ref{CIM_y1}) is compactly expressed in an operator format as
	
	\begin{equation} \label{CIM_ope_form}
		T{E_1} = \left( {{g_1}{\mathcal{Y}_1} - {g_{1n'}}} \right){E_1},
	\end{equation}
	where ${g_1}$, ${g_{1n'}}$ denote the integral operators for the Green’s function and its normal derivative on ${\gamma _1}$, respectively. Their expressions are given by
	\begin{align} \label{g1_gn1}
		{g_1} &= \oint_{{\gamma _1}} {{G_1}\left( {\mathbf{r},\mathbf{r}'} \right)} \left( . \right)dr',\notag \\
		g_{1n'} &= \oint_{{\gamma _1}} {\left[ {\frac{{\partial {G_1}\left( {\mathbf{r},\mathbf{r}'} \right)}}{{\partial n'}}} \right]} \left( . \right)dr'.\notag
	\end{align}
	After moving the second term of (\ref{CIM_ope_form}) on its right hand side (RHS) to its left hand side (LHS) and inverting the integral operator ${g_1}$, ${\mathcal{Y}_1}$ can be expressed as
	
	\begin{equation} \label{y1_ope_form_o}
		{\mathcal{Y}_1}{\rm{ = }}g_1^{ - 1}\left( {TI + {g_{1n'}}} \right),
	\end{equation}
	where $I$ is the identity operator. Therefore, (\ref{TH1_o}) can be rewritten via the SAO as
	\begin{equation} \label{TH1_o_y1}
		{{H_{t_1}}\left( {\mathbf{r}} \right)} = \frac{1}{{j\omega {\mu _1}}}{\mathcal{Y}_1}{\left. {{E_1}} \right|_{\mathbf{r} \in {\gamma _1}}}.
	\end{equation}
	
	Then, the surface equivalence theorem is applied through replacing the penetrable object inside ${\gamma_1}$ by its surrounding medium. The equivalent model with only the single electric current density is derived, as shown in Fig. \ref{single_model}(b), and $\widehat{\mathcal{Y}}_1$ can be expressed as
	\begin{equation} \label{y1_ope_form_e}
		{\widehat{\mathcal{Y}}_1}{\rm{ = }}{\widehat g}_1^{ - 1}\left( {TI + {{\widehat g}_{1n'}}} \right).
	\end{equation}
	The tangential magnetic field ${{\widehat{H}_{t_1}}\left( {\mathbf{r}} \right)}$ on ${\gamma_1}$ can also be achieved through the SAO in the equivalent model
	
	\begin{equation} \label{TH1_e_y1}
		{{{\widehat H}_{t_1}}\left( {\mathbf{r}} \right)} = \frac{1}{{j\omega {\mu _0}}}{\widehat{\mathcal{Y}}_1}{\left. {{{\widehat E}_1}} \right|_{\mathbf{r} \in {\gamma _1}}}.
	\end{equation}
	
	\begin{figure}
		\centering
		\begin{minipage}[h]{0.48\linewidth}
			\centerline{\includegraphics[scale=0.17]{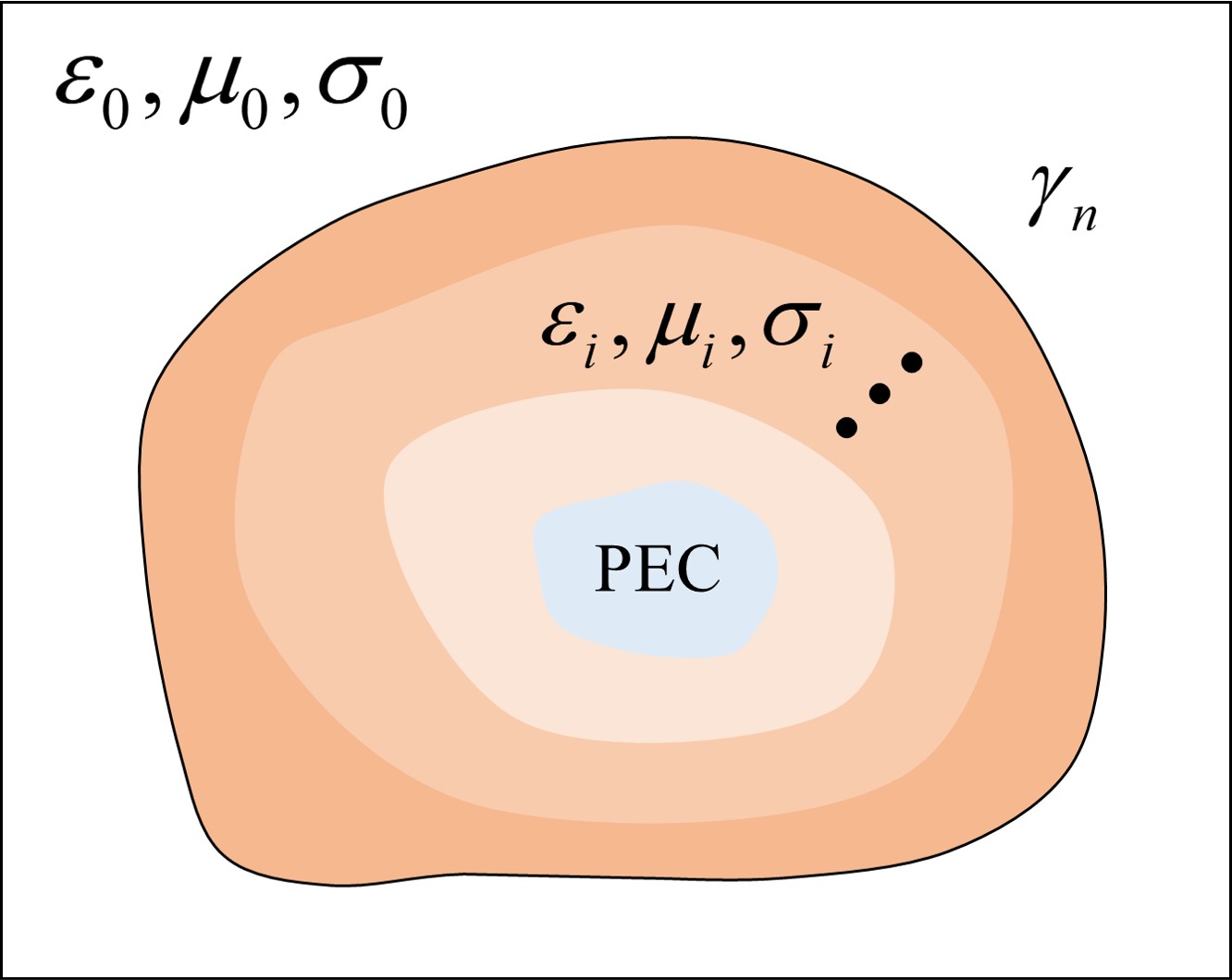}}
			\centerline{(a)}
		\end{minipage}
		\hfill
		\begin{minipage}[h]{0.48\linewidth}
			\centerline{\includegraphics[scale=0.17]{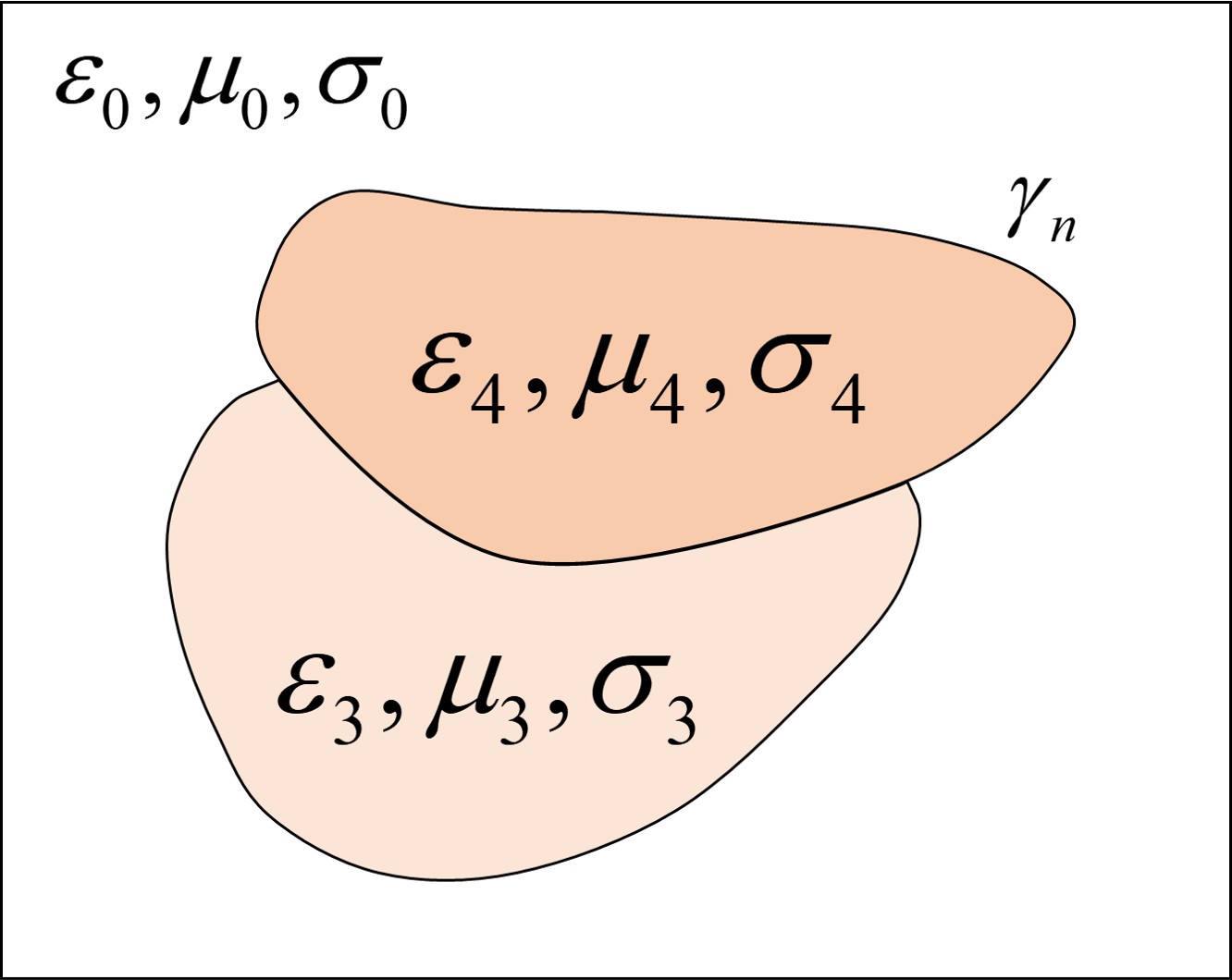}}
			\centerline{(b)}
		\end{minipage}
	\\	
	
		\begin{minipage}[h]{0.48\linewidth}
			\centerline{\includegraphics[scale=0.17]{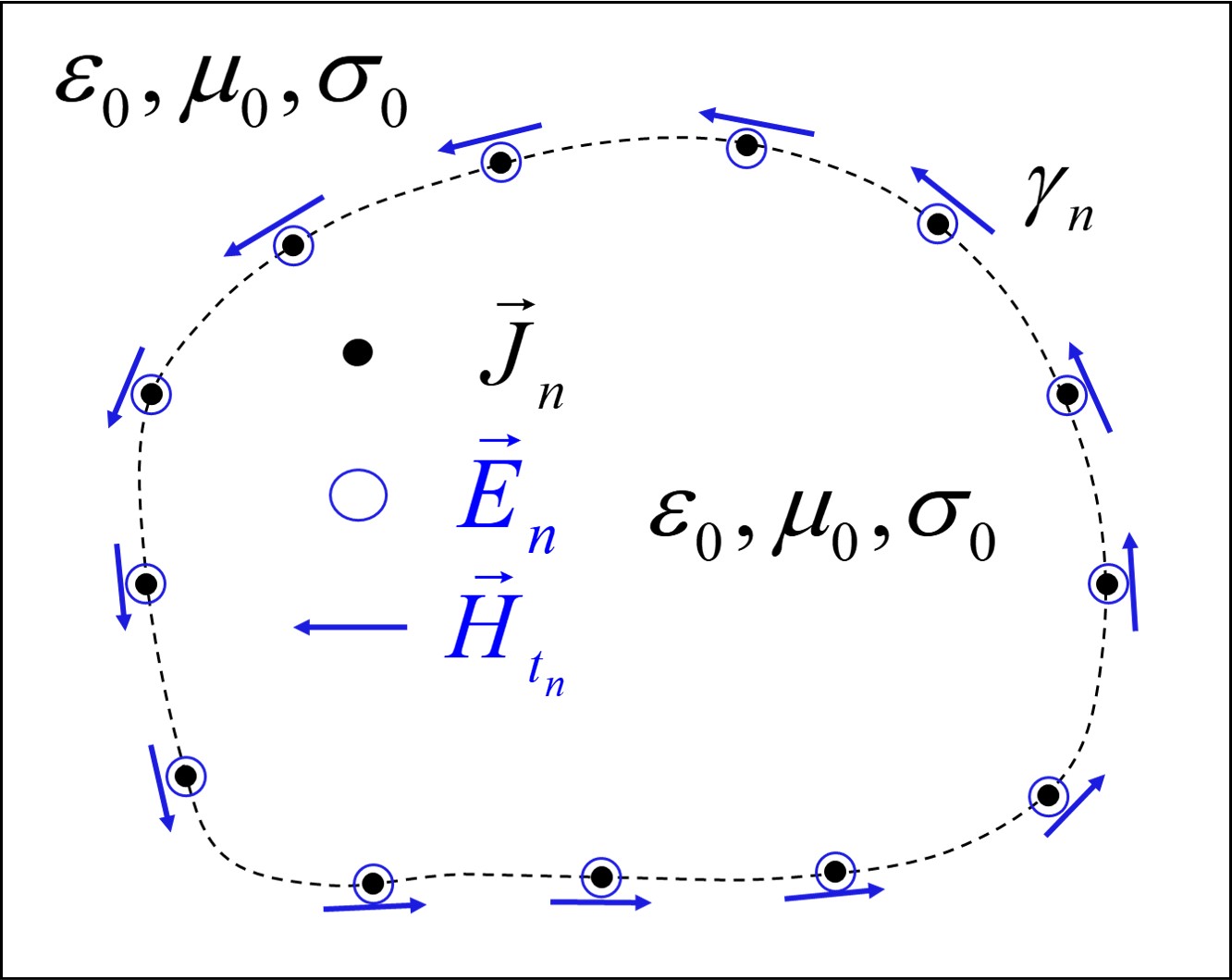}}
			\centerline{(c)}
		\end{minipage}
		\hfill
		\begin{minipage}[h]{0.48\linewidth}
			\centerline{\includegraphics[scale=0.17]{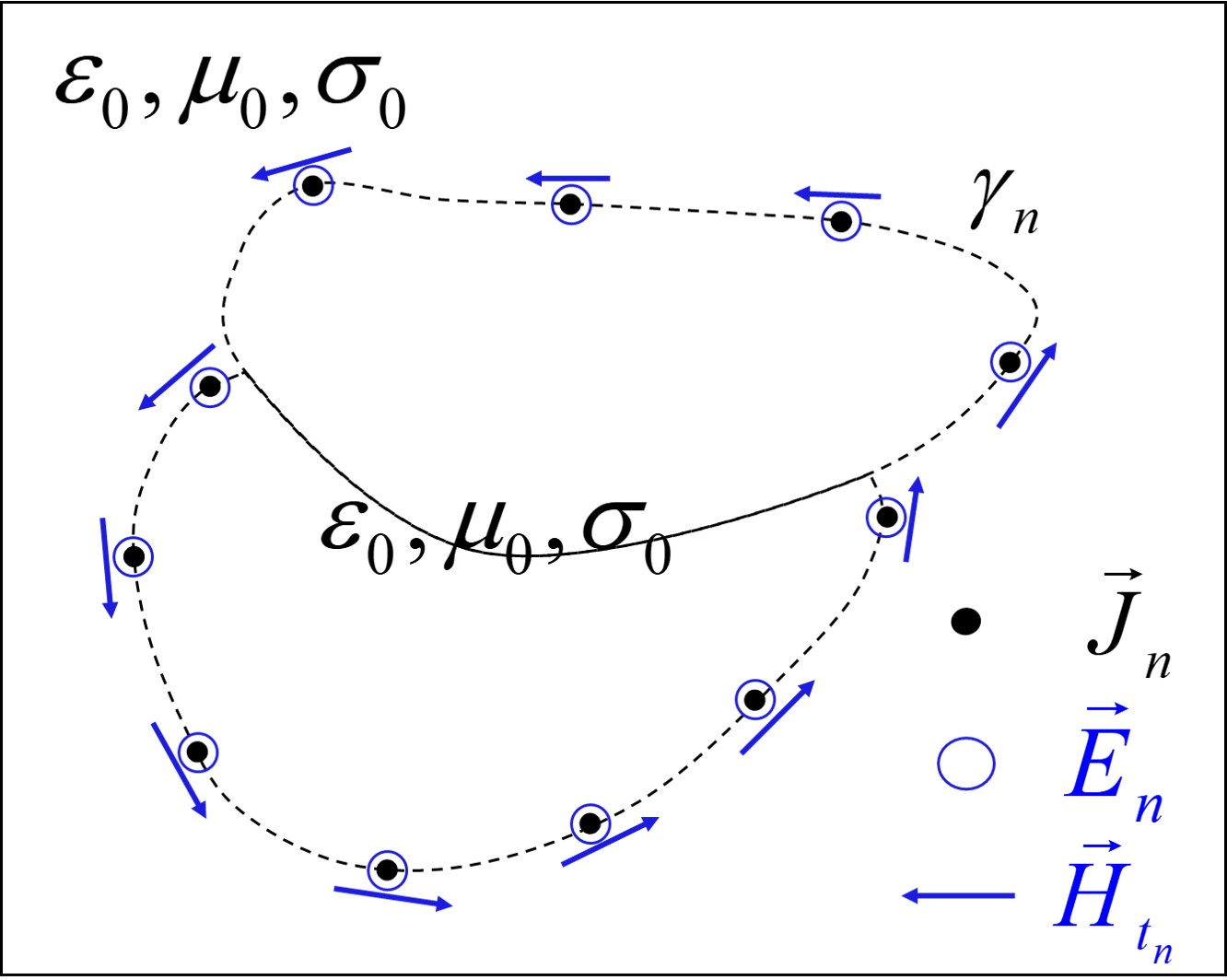}}
			\centerline{(d)}
		\end{minipage}
		\caption{(a) The original object embedded into multilayer dielectric media, (b) the original partial connected penetrable objects, (c) the equivalent model for object embedded into multilayer dielectric media, and (d) the equivalent model for partial connected penetrable objects.}
		\label{multi_partial_model} 
	\end{figure}
	
	Since the boundary condition (\ref{BCE1}) is enforced, the magnetic current density in the equivalent configuration vanishes, and only the electric current density ${J_1}$ is introduced on ${\gamma_1}$ to keep fields outside ${\gamma_1}$ unchanged. According to the surface equivalence theorem [\citen{LOVE}, Ch. 12, pp. 653-658], the equivalent electric current density ${J_1}$ is expressed as
	
	\begin{equation} \label{J1_H1}
		{J_1} = {{H_{t_1}}\left( {\mathbf{r}} \right)} - {{{\widehat H}_{t_1}}\left( {\mathbf{r}} \right)}.
	\end{equation}
	
	By substituting (\ref{TH1_o_y1}) and (\ref{TH1_e_y1}) into (\ref{J1_H1}), ${J_1}$ can be rewritten as
	\begin{equation} \label{J1_ys1}
		{J_1} = \underbrace {\left( {\frac{1}{{j\omega {\mu _1}}}{\mathcal{Y}}_1} - \frac{1}{{j\omega {\mu _0}}}{{\widehat{\mathcal{Y}}_1}} \right)}_{{\mathcal{Y}_{s_1}}}{\left. {{E_1}} \right|_{\mathbf{r} \in {\gamma _1}}},
	\end{equation}
	where ${\mathcal{Y}_{{s_1}}}$ is the differential surface admittance operator (DSAO) [\citen{Patel2016contour}] defined on $\gamma_1$. (\ref{J1_ys1}) relates the surface equivalent electric current density ${J_1}$ to electric fields ${E_1}$ on ${\gamma_1}$.
	
	When complex objects embedded in multilayers, composite structures with partial connected penetrable and PEC objects are involved in the SIE domain, as shown in Fig. \ref{multi_partial_model}(a) and (b), equivalent models with only the electric current density ${J_n}$ can still be derived through replacing those objects by the background medium, as shown in Fig. \ref{multi_partial_model}(c) and (d). Interested readers are referred to [\citen{Zhou2021embedded}][\citen{Zhou2021SS-SIE}] for more details to obtain the equivalent models. 
	
	In a word, the tangential magnetic fields on ${\gamma_n}$ for the original and equivalent models can also be achieved through the SAO given by
	
	\begin{equation} \label{THn_o_yn}
		{{H_{t_n}}\left( {\mathbf{r}} \right)}= \frac{1}{{j\omega {\mu _n}}}{\mathcal{Y}_n}{\left. {{E_n}} \right|_{\mathbf{r} \in {\gamma _n}}},
	\end{equation}
	
	\begin{equation} \label{THn_e_yn}
		{{{\widehat{H}_{t_n}}\left( {\mathbf{r}} \right)}} = \frac{1}{{j\omega {\mu _0}}}{\widehat {\mathcal{Y}}_n}{\left. {{{\widehat E}_n}} \right|_{\mathbf{r} \in {\gamma _n}}}.
	\end{equation}
	Here ${H_{t_n}}$, ${\widehat{H}_{t_n}}$, ${E_n}$, ${{\widehat E}_n}$, $\mathcal{Y}_n$, $\widehat {\mathcal{Y}}_n$ denote the tangential magnetic fields and electric fields on ${\gamma _n}$, and the SAO for the original and equivalent models, respectively.    
	
	Similar to the boundary condition (\ref{BCE1}), we also enforce ${\left. {{E_n}\left( {\mathbf{r}} \right)} \right|_{\mathbf{r} \in {\gamma _n}}} = {\left. {{{\widehat E}_n}\left( {\mathbf{r}} \right)} \right|_{\mathbf{r} \in {\gamma _n}}}$ on ${\gamma _n}$. Then, the equivalent electric current density on ${\gamma _n}$ can be derived and expressed as
	
	\begin{equation} \label{Jn_yn}
		{J_n} = \underbrace {\left( {\frac{1}{{j\omega {\mu _n}}}{\mathcal{Y}_n} - \frac{1}{{j\omega {\mu _0}}}{\widehat{\mathcal{Y}}_{n}}} \right)}_{\mathcal{Y}_{s_n}}{\left. {{E_n}} \right|_{\mathbf{r} \in {\gamma _n}}},
	\end{equation}
	where ${\mathcal{Y}}_{s_n}$ is the DSAO defined on $\gamma_n$.
	
	\subsection{The PDE Formulation in the PDE Domain}
	In the PDE domain, inhomogeneous or non-isotropic media are possibly included. The electric fields ${E_0}$ satisfy the scalar inhomogeneous Helmholtz equation
	\begin{equation} \label{IhoHe}
		{\nabla ^2}{E_0}{\rm{ + }}{k^2}{E_0} = j\omega \mu J,
	\end{equation}
	where $J$ denotes the electric current density inside the computational domain, and $\mu$ and $k$ denote the permeability and the wave number of the media inside the PDE domain, which may be a function of positions.
	
	(\ref{IhoHe}) can be expanded into its component form in TM mode as
	\begin{equation} \label{x_y_expand_Iho}
		\frac{\partial }{{\partial x}}\left( {\frac{1}{{{\mu _r}}}\frac{{\partial {E_0}}}{{\partial x}}} \right){\rm{ + }}\frac{\partial }{{\partial y}}\left( {\frac{1}{{{\mu _r}}}\frac{{\partial {E_0}}}{{\partial y}}} \right){\rm{ + }}k_0^2{\varepsilon _r}{E_0} = j\omega {\mu _0}J,
	\end{equation}
	where ${k_0}$ is the wave number of the free space, and ${\varepsilon _r}$, ${\mu _r}$  are the relative permittivity and permeability of the media inside the PDE domain, respectively.
	
	\subsection{The Hybrid SIE-PDE Formulation for the Whole Computational Domain}
	The hybrid SIE-PDE formulation can be derived through the following procedures.
	
	\begin{enumerate}
		\item In the SIE domain, the surface equivalence theorem is used to derive an equivalent model with only the electric current density ${J_n}$ on ${\gamma_n}$ by incorporating the DSAO. It is general and applicable for the single penetrable objects, objects embedded in multilayers, and arbitrarily connected penetrable and PEC objects [\citen{Zhou2021embedded}][\citen{Zhou2021SS-SIE}].  
		
		\item After that, the electric current density ${J_n}$ is used as an excitation in the inhomogeneous Helmholtz equation to couple the SIE and PDE formulations. Therefore, (\ref{x_y_expand_Iho}) is modified as
		
		\begin{align} 
			\frac{\partial }{{\partial x}}\left( {\frac{1}{{{\mu _r}}}\frac{{\partial {E_0}}}{{\partial x}}} \right){\rm{ + }}\frac{\partial }{{\partial y}}\left( {\frac{1}{{{\mu _r}}}\frac{{\partial {E_0}}}{{\partial y}}} \right) 
			{\rm{ + }}k_0^2{\varepsilon _r}{E_0} \label{Jn_modified} - j\omega {\mu _0}{J_n} = j\omega \mu J.
		\end{align}
		where $\mu_0$ is the permeability of the background medium. By substituting (\ref{Jn_yn}) into (\ref{Jn_modified}), it can be further rewritten as
		
		\begin{align} 
			\frac{\partial }{{\partial x}}\left( {\frac{1}{{{\mu _r}}}\frac{{\partial {E_0}}}{{\partial x}}} \right){\rm{ + }}\frac{\partial }{{\partial y}}&\left( {\frac{1}{{{\mu _r}}}\frac{{\partial {E_0}}}{{\partial y}}} \right)  
			\label{Jn+ysn_modified}{\rm{ + }}k_0^2{\varepsilon _r}{E_0}- j\omega {\mathcal{Y}_{s_n}}{\left.{{E_n}} \right|_{{\mathbf{r}} \in {\gamma _n}}} = j\omega \mu J.
		\end{align}

	\end{enumerate}
	
	Therefore, the governing equation for the whole computational domain is finally derived. It can be found that the SIE and PDE formulations are coupled through the surface equivalent electric current density, which does not require additional boundary conditions, e.g., the transmission condition, like other techniques, such as the FEM-MoM formulation [\citen{IlicFEM-MoM_ante_2009}-\citen{JiFEM/MoM2002}], the FE-BI formulation [\citen{EibertFE/BI1999}-\citen{YangFE-BI-MLFMA2013}], and the multi-region multi-solver [\citen{Guanmultisolver2017}][\citen{Guanmultisolver2016}]. The proposed hybrid SIE-PDE formulation is mathematically equivalent to the original physical model. Meanwhile, by coupling the SIE formulation, relatively coarse meshes are allowed when dealing with computationally challenging structures and less unknowns compared with the traditional FEM are required, which only reside on the boundary of the SIE domain. Therefore, through solving (\ref{Jn+ysn_modified}), all the electric fields of the computational domain can be accurately calculated. It should be noted that although ${E_0}$ exists in the whole computational domain, they are invalid in the homogeneous background medium replaced domain, in which fictitious fields are obtained, since the surface equivalence theorem is used. In the later section, we will discuss how to recover the electric fields inside the SIE domain. 
	
	\section{Detailed Implementations of the Proposed Hybrid SIE-PDE Formulation}
	\subsection{Construction of the Basis Functions for Electric Field Expansion}
	Before we discretize (\ref{Jn+ysn_modified}), the basis functions should be carefully selected to expand the electric fields. As is illustrated in Section III-C, the computational domain is decomposed into two \emph{overlapping} domains, and ${\gamma_n}$ is the boundary of the SIE domain, where ${J_n}$ resides in Fig. \ref{multi_partial_model}(c) and (d). In (\ref{Jn_yn}), we can find that ${\mathcal{Y}_{s_n}}$ is used to relate the surface equivalent current density ${J_n}$ to the tangential electric fields ${E_n}$ on ${\gamma_n}$. In our implementation, the Galerkin scheme is used to solve (\ref{Jn+ysn_modified}). The challenging part is that the discretized hybrid SIE-PDE formulation should be compatible on the boundary ${\gamma_n}$, so that the electromagnetic effects of ${J_n}$ can be directly imposed on the discretized Helmholtz equation. Since the linear basis functions are used to expand the electric fields in the PDE domain, the rooftop basis functions, which is constructed through combining the degenerated linear basis functions on ${\gamma_n}$ from two adjacent triangle elements, is used to expand the equivalent current density ${J_n}$ in Fig. \ref{multi_partial_model}(c) and (d). 
	
	The linear basis function in the two-dimensional space is defined on the triangular elements [\citen{Jin2015FEM}, Ch. 4, pp. 162-163] as
	\begin{equation} \label{linear_fun2}
		\begin{array}{*{20}{c}}
			{N_i^e\left( {x,y} \right) = \displaystyle \frac{1}{{2{\Delta ^e}}}\left( {a_j^e + b_j^ex + c_j^ey} \right),}&{i = 1,2,3},
		\end{array}
	\end{equation}
	where
	\[\begin{array}{*{20}{l}}
		{a_1^e = x_2^ey_3^e - y_2^ex_3^e}&{b_1^e = y_2^e - y_3^e}&{c_1^e = x_3^e - x_2^e},\\
		{a_2^e = x_3^ey_1^e - y_3^ex_1^e}&{b_2^e = y_3^e - y_1^e}&{c_2^e = x_1^e - x_3^e},\\
		{a_3^e = x_1^ey_2^e - y_1^ex_2^e}&{b_3^e = y_1^e - y_2^e}&{c_3^e = x_2^e - x_1^e},
	\end{array}\]
	and
	\[{\Delta ^e} = \frac{1}{2}\left| {\begin{array}{*{20}{c}}
			1&{x_1^e}&{y_1^e}\\
			1&{x_2^e}&{y_2^e}\\
			1&{x_3^e}&{y_3^e}
	\end{array}} \right| = \frac{1}{2}\left( {b_1^ec_2^e - b_2^ec_1^e} \right),\]
	where $x_i^e$ and $y_i^e$ ($i = 1,2,3$) denote the coordinates in the $x$, $y$ direction of the $i$th node of the $e$th triangle element, respectively, and ${\Delta ^e}$ denotes the area of the $e$th element. It should be noted that $N_i^e$ is locally supported, and only defined on the $e$th triangle element. On the boundary segment, the linear basis functions are degenerated into the following shape functions
	\begin{equation} \label{linear_fun1}
		\begin{array}{*{20}{c}}
			{N_1^{{l_n}}{\rm{ = }}\displaystyle \frac{{\left| {{\mathbf{r}_{n + 1}} - \mathbf{r}} \right|}}{{{l_n}}},}&{N_2^{{l_n}}{\rm{ = }}\displaystyle \frac{{\left| {\mathbf{r} - {\mathbf{r}_n}} \right|}}{{{l_n}}},}&{\mathbf{r} \in \left[ {{\mathbf{r}_n},{\mathbf{r}_{n + 1}}} \right]}. \notag
		\end{array}
	\end{equation}
	
	\begin{figure}
		\centering
		\includegraphics[width=0.45\textwidth]{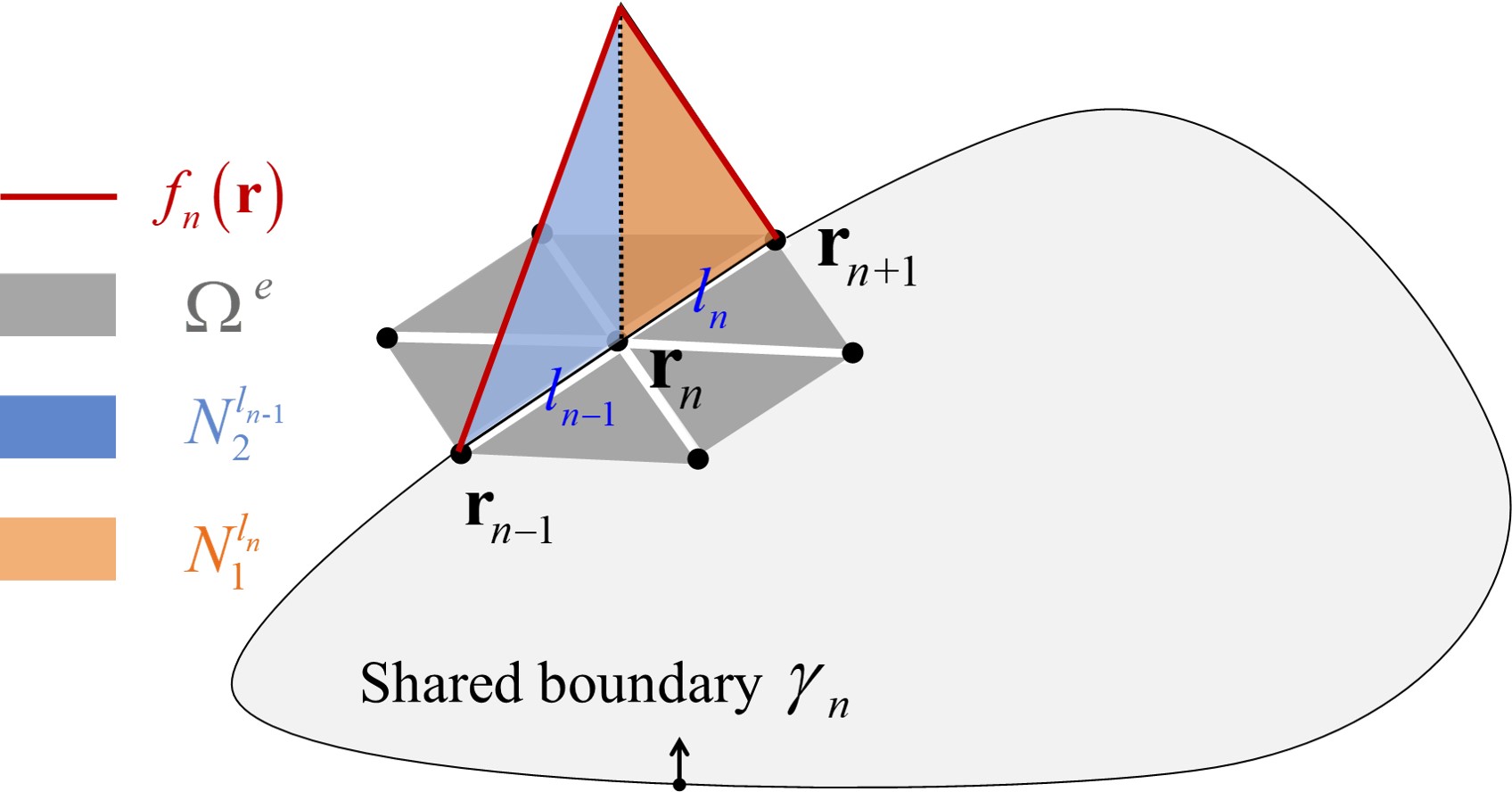}
		\caption{A general configuration to demonstrate the relationship between the linear basis function and the rooftop basis function on the share boundary $\gamma_n$. $\mathbf{r}_{n-1}$, $\mathbf{r}_n$, $\mathbf{r}_{n+1}$ are the endpoints of the $n$th and $(n + 1)$th segments, ${l_{n - 1}} = \left| {{\mathbf{r}_n} - {\mathbf{r}_{n - 1}}} \right|$, ${l_n} = \left| {{\mathbf{r}_{n + 1}} - {\mathbf{r}_n}} \right|$.}
		\label{rooftop} 
	\end{figure}
	
	To make the discretized SIE formulation and the PDE formulation compatible on $\gamma_n$, the rooftop basis function through combining two adjacent triangle elements is defined at the node rather than on the segment, which can be expressed at the $n$th node as
	\begin{equation} \label{roogtop_fun}
		{f_n}\left( {\bf{r}} \right) = \left\{ {\begin{array}{*{20}{c}}
				{\displaystyle \frac{{\left| {\mathbf{r} - {\mathbf{r}_{n - 1}}} \right|}}{{{l_n}}},}&{{\mathbf{r}} \in \left[ {{\mathbf{r}_{n - 1}},{\mathbf{r}_n}} \right]}\\
				{\displaystyle \frac{{\left| {{\mathbf{r}_{n + 1}} - \mathbf{r}} \right|}}{{{l_n}}},}&{{\mathbf{r}} \in \left[ {{\mathbf{r}_n},{\mathbf{r}_{n + 1}}} \right]}
		\end{array}} \right..
	\end{equation}
	Therefore, as shown in Fig. \ref{rooftop}, the relationship between the linear basis function and the rooftop basis function can be expressed as
	
	\begin{equation} \label{linear-rooftop-relation} 	
		{f_n}\left( {\bf{r}} \right) = \left\{ {\begin{array}{*{20}{l}}
				{N_2^{{l_{n - 1}}},}&{{\mathbf{r}} \in \left[ {{\mathbf{r}_{n - 1}},{\mathbf{r}_n}} \right]}\\
				{N_1^{{l_n}},}&{{\mathbf{r}} \in \left[ {{\mathbf{r}_n},{\mathbf{r}_{n + 1}}} \right]}
		\end{array}} \right..
	\end{equation} 
	
	\subsection{Discretization for the Proposed Hybrid SIE-PDE Formulation}
	To clearly explain the procedure of discretization for the proposed hybrid SIE-PDE formulation, let us consider the general scenario with a single penetrable object in Fig. \ref{single_model}(a). The penetrable object with the boundary $\gamma_1$ is located in the SIE domain, and its equivalent model is illustrated in detail in Section III-A.
	
	We first discretize $\gamma_1$ into $m_1$ segments and use the rooftop basis function to expand $E_1$ and $H_{t_1}$ in (\ref{CIM}) and (\ref{TH1_o}) as
	
	\begin{equation} \label{E1_fn}
		{E_1}\left( {\mathbf{r}} \right) = \sum\limits_{n = 1}^{{m_1}} {{e_{1n}}{f_n}\left( {\mathbf{r}} \right)}, 
	\end{equation}
	
	\begin{equation} \label{H1_fn}
		{{H_{t_1}}\left( {\mathbf{r}} \right)} = \sum\limits_{n = 1}^{m_1} {{h_{1n}}{f_n}\left( {\mathbf{r}} \right)},
	\end{equation}
	where ${f_n}\left( {\mathbf{r}} \right)$ denotes the $n$th basis function defined in (\ref{roogtop_fun}). After substituting (\ref{E1_fn}) and (\ref{H1_fn}) into (\ref{CIM}) and testing (\ref{CIM}) through the Galerkin scheme at each segment node, we collect all $E_1$ and $H_{t_1}$ expansion coefficients into two column vectors ${\mathbf{E}_1}$ and ${\mathbf{H}_1}$ as
	\begin{equation} \label{E1_column}
		{\mathbf{E}_1} = {\left[ {\begin{array}{*{20}{c}}
					{{e_{11}}}&{{e_{12}}}& \ldots &{{e_{1{m_1}}}}
			\end{array}} \right]^T},
	\end{equation}
	
	\begin{equation} \label{H1_column}
		{\mathbf{H}_1} = {\left[ {\begin{array}{*{20}{c}}
					{{h_{11}}}&{{h_{12}}}& \ldots &{{h_{1{m_1}}}}
			\end{array}} \right]^T}.
	\end{equation}
	Then, (\ref{CIM}) can be rewritten into the matrix form as
	
	\begin{equation} \label{CIM_matrix}
		T{\mathbb{L}_1}{\mathbf{E}_1} = {\mathbb{P}_1}{\mathbf{H}_1} +{\mathbb{U}_1}{\mathbf{E}_1},
	\end{equation}
	where the subscript 1 denotes that the testing procedure is applied on $\gamma_1$, and the entities of ${\mathbb{L}_1}$, ${\mathbb{P}_1}$ and ${\mathbb{U}_1}$ are expressed as
	\begin{align} 	
	&\label{L1} {\left[ {\mathbb{L}_1} \right]_{m,n}}= \int_{s^{m}} {{f_m}\left( {\mathbf{r}} \right)\int_{s^{n}} {{f_n}\left( {\mathbf{r}'} \right)dr'dr} },\\
	&\label{P1}	{\left[ {{\mathbb{P}_1}} \right]_{m,n}}= j\omega {\mu _1} \int_{s^{m}}{{f_m}\left( {\mathbf{r}} \right)\int_{s^{n}}\!{{G_1}\left( {\mathbf{r},\mathbf{r}'} \right){f_n}\left( {\mathbf{r}'} \right)dr'dr} },\\ 
	&\label{U1}	{\left[ {{\mathbb{U}_1}} \right]_{m,n}} =\!\int_{s^{m}}\! {{f_m}\left( {\mathbf{r}} \right)\!\int_{s^{n}}\!{{k_1}\frac{{\vec \rho  \cdot \hat n'}}{\rho }{G_1}^\prime\! \left( {\mathbf{r},\mathbf{r}'} \right)\!{f_n}\!\left( {\mathbf{r}'} \right)\!dr'\!dr} }, 
	\end{align}
	where ${G_1}^\prime  =  - jH_1^{\left( 2 \right)}\left( {{k_1}\rho } \right)/4$, ${s^m} = \gamma _1^{m - 1} + \gamma _1^m$, and ${s^n} = \gamma _1^{n - 1} + \gamma _1^n$. 
	\vspace{4 pt}
	
	In this paper, a modified singularity handling approach through combining the numerical Gaussian quadrature and analytical integration approaches is used to accurately and efficiently handle the singular and nearly singular integrals in (\ref{P1}) and (\ref{U1}). In this approach, the Hankel function and its gradient are evaluated by the small variable approximation when a predefined threshold is satisfied. Then, the analytical approach is used to compute the singular and nearly singular integration, and the Gaussian quadrature is applied in other non-singular integration. To make our paper more focused, the detailed techniques are listed in the appendix.
	
	Then, through moving the second term of (\ref{CIM_matrix}) on the RHS to its LHS and inverting the square coefficient matrix ${\mathbb{P}_1}$, we obtain the SAO ${\mathbb{Y}_1}$  [\citen{Patel2016contour}] as
	
	\begin{equation} \label{H1_Y1}
		{\mathbf{H}_1} = \underbrace {{{\left[ {\mathbb{P}_1} \right]}^{ - 1}}\left( {T{\mathbb{L}_1} - {\mathbb{U}_1}} \right)}_{{\mathbb{Y}_1}}{{\mathbf{E}_1}}.
	\end{equation}
	In the equivalent configurations in Fig. \ref{single_model}(b), ${\widehat{H}_{t_1}\left( {\mathbf{r}} \right)} $ can be also expanded through the rooftop basis function as
	
	\begin{equation} \label{H1_e_fn}
		{\widehat{H}_{t_1} \left( {\mathbf{r}} \right)} = \sum\limits_{n = 1}^{{m_1}} {\widehat {h}_{1n}{f_n}\left( {\mathbf{r}} \right)}.
	\end{equation}
	With the similar procedures in the original object, we have
	
	\begin{equation} \label{H1_e_Y1}
		{\widehat{\mathbf{H}}_1} = \underbrace {{{\left[ {\widehat{\mathbb{P}}_1} \right]}^{ - 1}}\left( {T{\mathbb{L}_1} - {\widehat{\mathbb{U}}_1}} \right)}_{\widehat{\mathbb{Y}}{_1}}{{\mathbf{E}_1}},
	\end{equation}
	where
	\begin{equation} \label{H1_e_column}
		{\widehat{\mathbf{H}}_1} = {\left[ {\begin{array}{*{20}{c}}
					{{\widehat{h}_{11}}}&{{\widehat{h}_{12}}}& \ldots &{{\widehat{h}_{1{m_1}}}}
			\end{array}} \right]^T}.
	\end{equation}
	
	Since (\ref{BCE1}) is enforced on $\gamma_1$, only the surface equivalent electric current density ${J}_1$ is required on $\gamma_1$ in the equivalent configuration to keep fields outside $\gamma_1$ unchanged. It can also be expanded through the rooftop basis function, and the coefficients are collected into the column vector ${\mathbf{J}_1}$ as
	
	\begin{equation} \label{J1_column}
		{\mathbf{J}_1} = {\left[ {\begin{array}{*{20}{c}}
					{{j_{11}}}&{{j_{12}}}& \ldots &{{j_{1{m_1}}}}
			\end{array}} \right]^T}.
	\end{equation}
	By substituting (\ref{H1_Y1}) and (\ref{H1_e_Y1}) into (\ref{J1_H1}), ${\mathbf{J}_1}$ is obtained as
	
	\begin{equation} \label{J1_Ys1}
		{\mathbf{J}_1} = {{\mathbb{Y}}_{s_1}}{\mathbf{E}_1},
	\end{equation}
	where ${{\mathbb{Y}}_{s_1}}$ is the DSAO [\citen{Patel2016contour}], and can be expressed as
	\begin{align} 
		{{\mathbb{Y}}_{s_1}} = {{\mathbb{Y}}_1} - {\widehat{\mathbb{Y}}_1}
		\label{Ys1} = {{\left[ {\mathbb{P}_1} \right]}^{ - 1}}\left( {T{\mathbb{L}_1} - {\mathbb{U}_1}} \right) - {{\left[ {\widehat{\mathbb{P}}_1} \right]}^{ - 1}}\left( {T{\mathbb{L}_1} - {\widehat{\mathbb{U}}_1}} \right).
	\end{align}
	
	As shown in Fig. \ref{multi_partial_model}(a) and (b), when complex objects embedded in multilayers or composite objects with partial connected penetrable objects and PEC objects are involved in the SIE domain, $\mathbf{J}_n$ on $\gamma_n$ can also be derived according to the surface equivalence theorem [\citen{LOVE}, Ch. 12, pp. 653-658] in Fig. \ref{multi_partial_model}(c) and (d), which can be expressed as

	\begin{equation} \label{Jn_Ysn}
		{\mathbf{J}_n} = {{\mathbb{Y}}_{s_n}}{\mathbf{E}_n}.
	\end{equation}
	Detailed formulations of ${\mathbb{Y}}_{s_n}$ on $\gamma_n$ can be found in [\citen{Zhou2021embedded}][\citen{Zhou2021SS-SIE}].  Interested readers are referred to them for more details.
	
	The remaining computational domain and homogenous background medium replaced domain are treated as the PDE domain, and the scalar inhomogeneous Helmholtz equation in terms of the electric fields $E_0$ is expressed in (\ref{Jn_modified}). We apply the Galerkin scheme to deriving the discretized formulation. By moving the excitation term of (\ref{Jn_modified}) on the RHS to its LHS, the residual error is defined as
	
	\begin{align} 
		\label{residual_error1}r = \frac{\partial }{{\partial x}}\left( {\frac{1}{{{\mu _r}}}\frac{{\partial {E_0}}}{{\partial x}}} \right){\rm{ + }}\frac{\partial }{{\partial y}}\left( {\frac{1}{{{\mu _r}}}\frac{{\partial {E_0}}}{{\partial y}}} \right)&{\rm{ + }}k_0^2{\varepsilon _r}{E_0}- j\omega {\mu _0}{J_n} - j\omega \mu J.
	\end{align}
	
	The PDE domain is discretized into $M$ triangle elements, and the triangle meshes are compatible with segments on $\gamma_n$, as shown in Fig. \ref{rooftop}. In each triangle element, $E_0$ can be approximated as
	
	\begin{equation} \label{E0_exp}
		\begin{array}{*{20}{c}}
			{{E_0}^e\left( {x,y} \right) = \sum\limits_{j = 1}^3 {N_j^e\left( {x,y} \right){E_0}_j^e,} }&{e = 1,2,...,M}.
		\end{array}
	\end{equation}
	The weighted residual on the $e$th triangle element can be expressed as 
	\begin{equation} \label{Ri}
		\begin{array}{*{20}{c}}
			R_i^e = \iint_{{\Omega ^e}}{N_i^er{d}x{d}y},&i = 1,2,3,
		\end{array}
	\end{equation}
	where ${\Omega ^e}$ is the $e$th element, and $N_i^e$ are the two-dimensional linear basis functions given by (\ref{linear_fun2}). By substituting (\ref{residual_error1}) into (\ref{Ri}) and applying the divergence theorem [\citen{Jin2015FEM}, Ch. 4, pp. 154-157], (\ref{Ri}) can be rewritten as
	\begin{align} 
	\label{Ri+r} R_i^e = &{ \iint_{\Omega ^e}}\Big( - \frac{1}{{\mu _r^e}}\frac{{\partial N_i^e}}{{\partial x}}\frac{{\partial {E_0}}}{{\partial x}} - \frac{1}{{\mu _r^e}}\frac{{\partial N_i^e}}{{\partial y}}\frac{{\partial {E_0}}}{{\partial y}}\Big. \notag \\ 
	&\Big. + k_0^2\varepsilon _r^eN_i^e{E_0} \Big){d}x{d}y- j\omega {\mu _0} \iint_{{\Omega ^e}}{N_i^e}{J_n}{d}x{d}y - j\omega \mu \iint_{\Omega ^e}{N_i^e}J{d}x{d}y + \oint_{{\Gamma ^e}}{N_i^e} {\bf{D}} \cdot {\mathbf{n}^e}{d}\Gamma.&
	\end{align}

In (\ref{Ri+r}), $\varepsilon _r^e$, $\mu _r^e$ denote the relative permittivity and permeability of the $e$th triangle element, ${\Gamma ^e}$ denotes the contour of ${\Omega ^e}$, ${\mathbf{n}^e}$ is the unit normal vector pointing outward ${\Gamma ^e}$ , and
	
	\begin{equation} \label{D}
		{\mathbf{D}} = \left( {\frac{1}{{{\mu _r}}}\frac{{\partial {E_0}}}{{\partial x}}\mathbf{x} + \frac{1}{{{\mu _r}}}\frac{{\partial {E_0}}}{{\partial y}}\mathbf{y}} \right).
	\end{equation}

	As is illustrated in Section III, $J_n$ only exists on $\gamma_n$. It can also be approximated with the linear basis function on $\gamma_n$ as		
	\begin{equation} \label{Jn_exp}
		\begin{array}{*{20}{c}}
			{J_n^s\left( {x,y} \right) = \sum\limits_{j = 1}^2 {N_j^{l_s}\left( {x,y} \right){J_n}_{j}^s,} }&{s = 1,2,...,{m_n}}.
		\end{array}
	\end{equation}
	Therefore, the integration terms of $J_n$ is reduced into one-dimensional space. By substituting (\ref{E0_exp}), (\ref{Jn_exp}) into (\ref{Ri+r}), the elemental equation is expressed as 
	\begin{align}
		R_i^e =& {\sum\limits_{j = 1}^3 \iint_{\Omega ^e}}\Big(  - \frac{1}{{\mu _r^e}}\frac{{\partial N_i^e}}{{\partial x}}\frac{{\partial N_j^e}}{{\partial x}} - \frac{1}{{\mu _r^e}}\frac{{\partial N_i^e}}{{\partial y}}\frac{{\partial N_j^e}}{{\partial y}}\Big. \notag \\
		\label{Ri_Jn_exp}&\Big. + k_0^2\varepsilon _r^eN_i^eN_j^e \Big){E_0}_j^e{d}x{d}y - j\omega {\mu _0}\sum\limits_{j = 1}^2 {\left( {\int_{{\gamma _n^s}} {N_i^{l_s}N_j^{l_s}} } \right)} {J_n}_{j}^s{d}x{d}y - j\omega \mu \iint_{\Omega ^e}{N_i^e}J{d}x{d}y + \oint_{{\Gamma ^e}}{N_i^e} {\bf{D}} \cdot {{\mathbf{n}^e}{d}\Gamma }.&
	\end{align}
	and it can be written in the matrix form as  
	\begin{equation} \label{SIE-PDE_ele}
		{\mathbf{R}^e} =  - {\mathbb{K}^e}{\mathbf{E}^e} - {\mathbb{B}^s}{\mathbf{J}_n^s} - {\mathbf{b}^e} + {\mathbf{g}^e},
	\end{equation}
	where
	\begin{align}\notag
		&{{{\mathbf{R}}^e} = {{\left[ {\begin{array}{*{20}{c}}
							{R_1^e}&{R_2^e}&{R_3^e}
					\end{array}} \right]}^T},}\\ \notag
		&{{{\mathbf{E}}^e} = {{\left[ {\begin{array}{*{20}{c}}
							{{E_0}_{1}^e}&{{E_0}_{2}^e}&{{E_0}_{3}^e}
					\end{array}} \right]}^T},}\\ \notag
		&{{\mathbf{J}}_n^s = {{\left[ {\begin{array}{*{20}{c}}
							{{J_n}_{1}^s}&{{J_n}_{2}^s}
					\end{array}} \right]}^T}.}  \notag
	\end{align}
	The entities of ${\mathbb{K}^e}$, ${\mathbb{B}^s}$, ${\mathbf{b}^e}$ and ${\mathbf{g}^e}$ are given by

\begin{align}
	 \label{K} & {\left[ \mathbb{K}^e\right]_{i,j} = \iint_{\Omega ^e}}\Big( - \frac{1}{{\mu _r^e}}\frac{{\partial N_i^e}}{{\partial x}}\frac{{\partial N_j^e}}{{\partial x}} - \frac{1}{{\mu _r^e}}\frac{{\partial N_i^e}}{{\partial y}}\frac{{\partial N_j^e}}{{\partial y}} \Big. \Big.+ k_0^2\varepsilon _r^eN_i^eN_j^e \Big) {d}x{d}y,{\quad i,j = 1,2,3,} \\
	 \label{B} &{{{\left[ \mathbb{B}^s \right]}_{i,j}} = j\omega {\mu _0}\int_{{\gamma ^s}} {N_i^{l_s}N_j^{l_s}} dxdy,}{\quad i,j = 1,2,}\\
	 \label{b} &{{{\left[ {\mathbf{b}^e} \right]}_i} = j\omega \mu \iint_{\Omega ^e}{N_i^e}J{d}x{d}y,}{\quad i = 1,2,3,}\\
	 \label{g} &{{{\left[ {\mathbf{g}^e} \right]}_i} = \oint_{{\Gamma ^e}} {N_i^e} {\bf{D}} \cdot {\mathbf{n}^e}{d}\Gamma,}{\quad i = 1,2,3.}
\end{align}
	Then, the final matrix equation can be obtained by assembling (\ref{SIE-PDE_ele}) for all elements, given by
	
	\begin{equation} \label{Ri_all}
		{\mathbf{R}} =  -{\mathbb{K}}{\mathbf{E}} - {\mathbb{B}}{\mathbf{J}_n} - {\mathbf{b}} + {\mathbf{g}}.
	\end{equation}
	By enforcing that the residual error vanishes on all the elements, (\ref{Ri_all}) can be written as
	
	\begin{equation} \label{Ri_for_0}
		{\mathbb{K}}{\mathbf{E}} + {\mathbb{B}}{\mathbf{J}_n} = {\mathbf{g}} - {\mathbf{b}}.
	\end{equation}
	By substituting (\ref{Jn_Ysn}) into (\ref{Ri_for_0}), it can be rewritten as
	
	\begin{equation} \label{pre-SIE-PDE}
		{\mathbb{K}}{\mathbf{E}} + {\mathbb{B}}{\mathbb{Y}}_{s_n}{\mathbf{E}_n}= {\mathbf{g}} - {\mathbf{b}}.
	\end{equation}
	It should be noted that ${\left. {{E_n}\left( {\mathbf{r}} \right)} \right|_{\mathbf{r} \in {\gamma _n}}} = {\left. {{E_0}\left( {\mathbf{r}} \right)} \right|_{\mathbf{r} \in {\gamma _n}}}$. Therefore, the second term of (\ref{pre-SIE-PDE}) on the LHS can be expanded through replacing $\mathbf{E}_n$ by $\mathbf{E}$, and the discretized SIE-PDE formulation can be obtained as  
	
	\begin{equation} \label{SIE-PDE}
		\left( {\mathbb{K}} + {\mathbb{A}} \right) {\mathbf{E}}= {\mathbf{g}} - {\mathbf{b}},
	\end{equation}
	where $\mathbb{A}$ denotes the expansion matrix for ${\mathbb{B}}{\mathbb{Y}}_{s_n}$. Through solving (\ref{SIE-PDE}), we can obtain the electric fields in the PDE domain. Other interested parameters, like near fields, the radar cross section (RCS) can be calculated. In the following subsection, we briefly discuss how to calculate the electric fields in the whole computational domain. 
	
	\subsection{Fields Calculation in the Whole Computational Domain}
	By solving (\ref{SIE-PDE}), the electric fields in the PDE domain are obtained. However, since the surface equivalence theorem is applied, the electric fields inside the SIE domain cannot be derived directly from the proposed SIE-PDE formulation. Therefore, electric fields inside the SIE domain should be carefully recovered. As illustrated in Section IV-B, if the objects are embedded in multilayers, or composite objects with partially connected penetrable objects and PEC objects are involved, electric fields inside the SIE domain can be calculated through the reverse procedure used to construct the equivalent model. More details can be found in [\citen{Zhou2021embedded}][\citen{Zhou2021SS-SIE}]. 

	\section{Numerical Results and Discussion}
	All numerical examples in this paper are performed on a workstation with an Intel i7-7700 3.6 GHz CPU and 64 G memory. Our in-house codes are written in Matlab. To make a fair comparison, only a single thread without any parallel computation is used to complete the simulations.

	\subsection{Accuracy Verification of the Proposed Hybrid SIE-PDE Formulation for Penetrable Structures}	
	To verify the accuracy of the proposed SIE-PDE formulation, we first consider TM scattering problems induced by simple penetrable objects. The homogeneous background medium is chosen to simplify the verification. Two dielectric objects, an infinitely long dielectric cylinder and cuboid with the relative permittivity ${\varepsilon _r} = 2.3$ are selected. Both objects are located in the air. The first-order absorbing boundary condition [\citen{Jin2015FEM}, Ch. 4, pp. 207-212] is used to truncate the computational domain in the FEM and the proposed SIE-PDE formulation, which is expressed as
	\begin{equation} \label{ABC}
		\frac{1}{{{\mu _r}}}\frac{{\partial {E_0}}}{{\partial n}} + \gamma {E_0} = q,
	\end{equation}
	with $\gamma$ and $q$  given by
	\[\gamma  = \frac{1}{{{\mu _r}}}\left[ {j{k_0} + \frac{{\kappa \left( s \right)}}{2}} \right],\]
	\[q = \frac{1}{{{\mu _r}}}\frac{{\partial {E^{inc}}}}{{\partial n}} + \frac{1}{{{\mu _r}}}\left[ {j{k_0} + \frac{{\kappa \left( s \right)}}{2}} \right]{E^{inc}},\]
	where $\kappa \left( s \right)$  is the curvature of the truncated boundary at the point $s$, $\mu_r$  is the relative permeability of the surrounding medium, and ${E^{inc}}$ is the incident plane wave. In our simulation, the truncated boundary is a circle with a radius of 6 m. It is 5 m away from the dielectric objects to reduce possible reflection from the truncated boundary, and $\kappa \left( s \right)$  is 1/6. 
		\begin{figure}
		\centering
		\begin{minipage}[t]{0.48\linewidth}
			\centering
			\centerline{\includegraphics[scale=0.34]{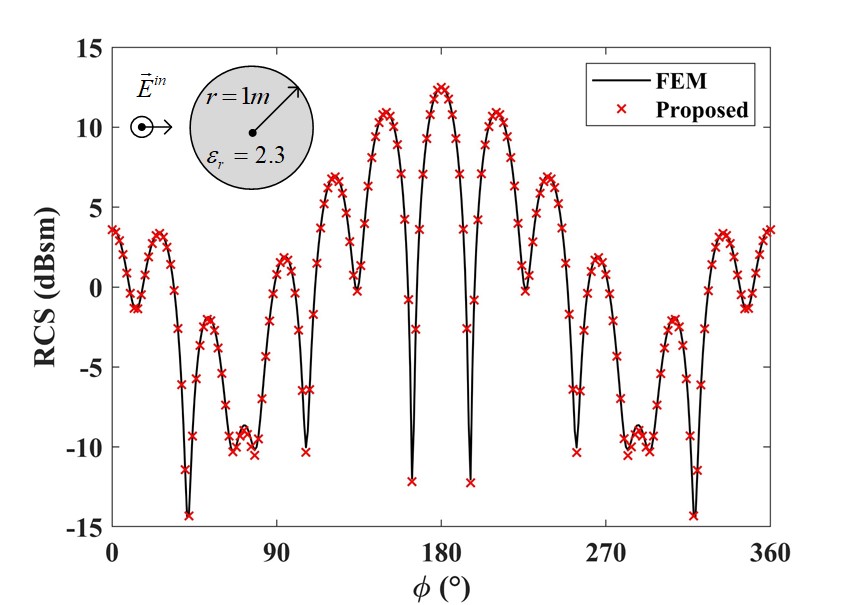}}
			\centering
			\centerline{(a)}
		\end{minipage}
		\hfill
		\begin{minipage}[t]{0.48\linewidth}
			\centering
			\centerline{\includegraphics[scale=0.34]{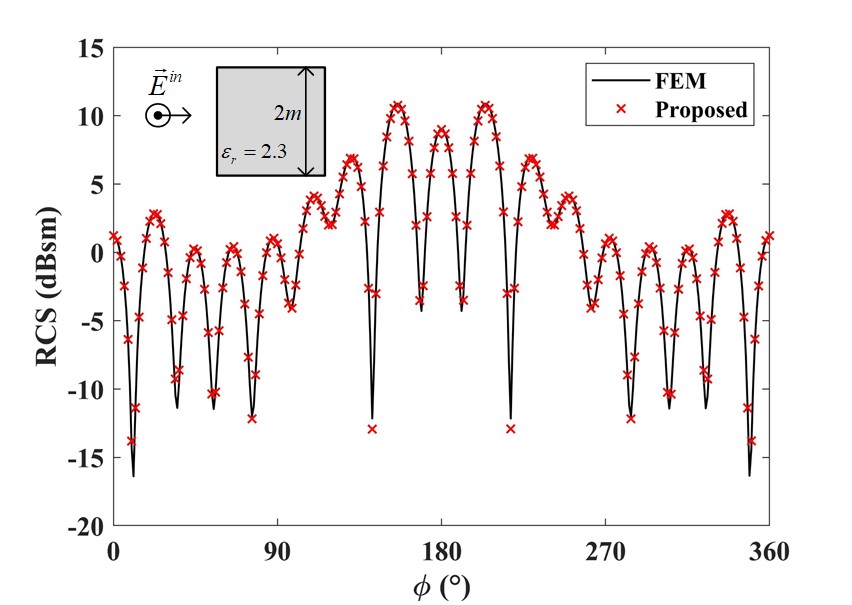}}
			\centering
			\centerline{(b)}
		\end{minipage}
		\caption{The RCS obtained from the proposed SIE-PDE formulation and the FEM for (a) the dielectric cylinder, and (b) the dielectric cuboid.}
		\label{RCS}
	\end{figure}
	
	\begin{figure}
		\centering
		\begin{minipage}[t]{0.48\linewidth}
			\centerline{\includegraphics[scale=0.062]{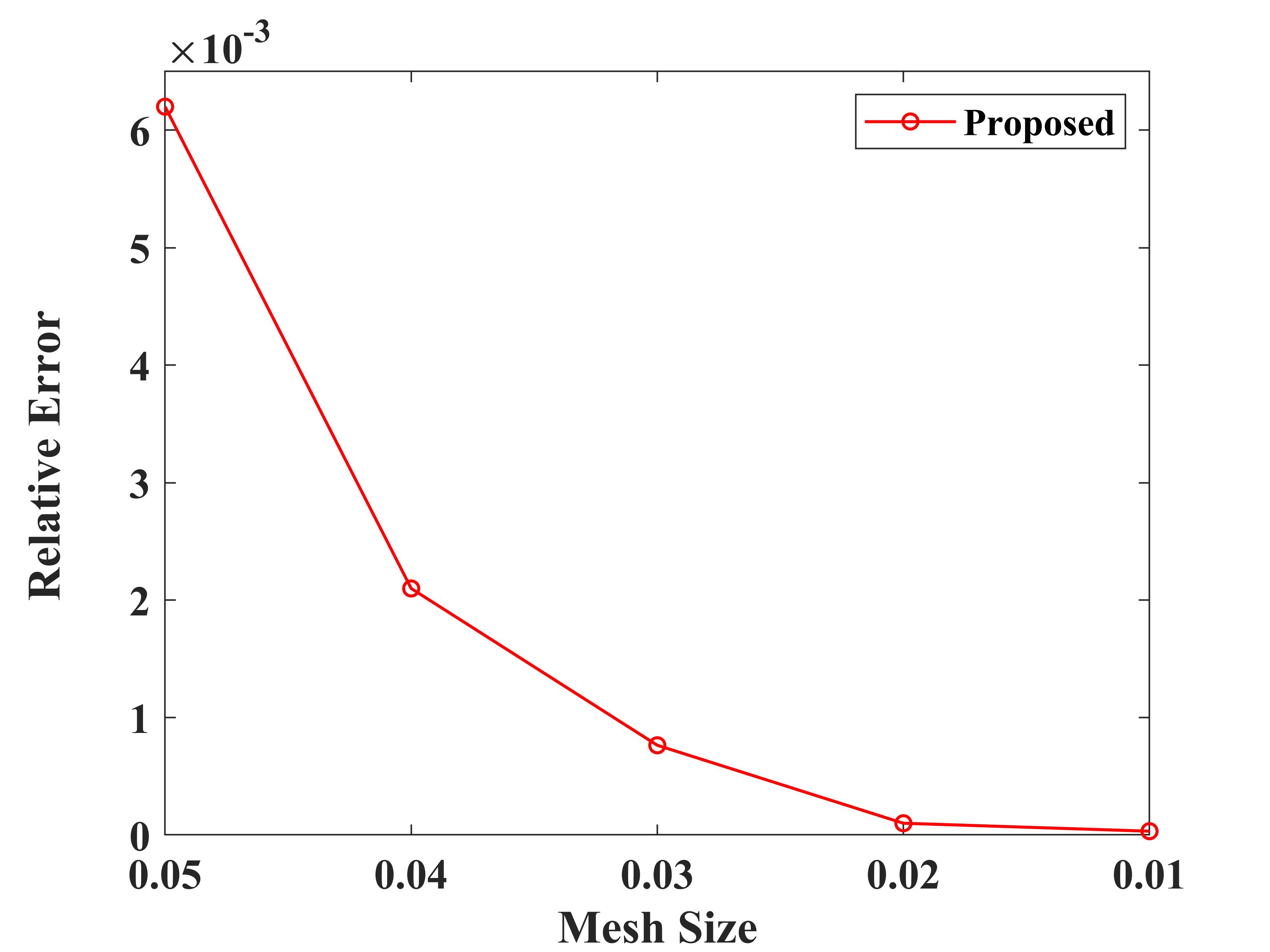}}
			\centerline{(a)}
		\end{minipage}
		\hfill
		\begin{minipage}[t]{0.48\linewidth}
			\centerline{\includegraphics[scale=0.062]{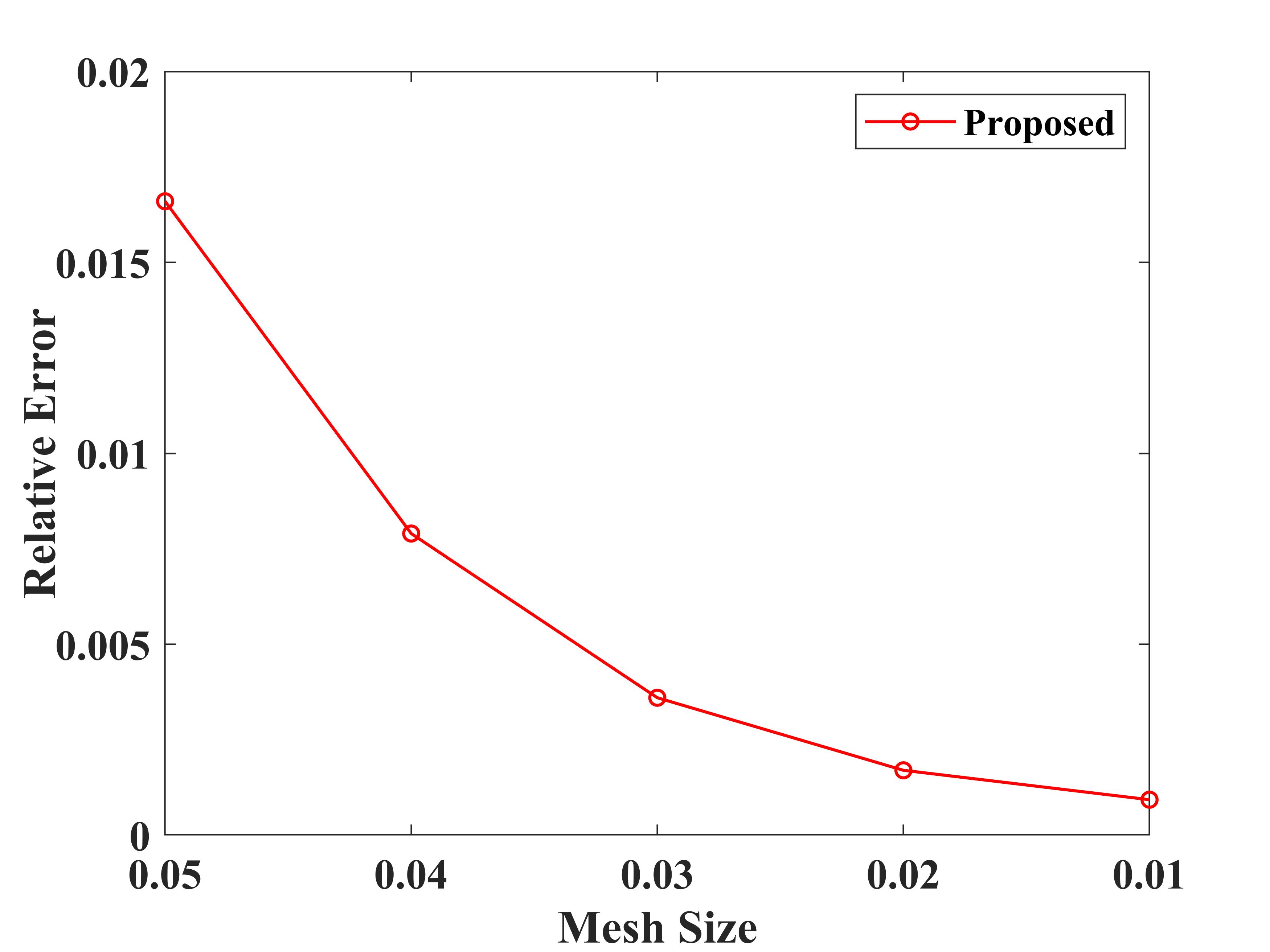}}
			\centerline{(b)}
		\end{minipage}
		\caption{Relative error of the proposed SIE-PDE formulation for (a) the dielectric cylinder, and (b) the dielectric cuboid.}
		\label{convergence}
	\end{figure}
	The radius of the cylinder is 1 m, and the side length of the cuboid is 2 m. A TM-polarized plane wave incidents from the $x$-axis with $f = 300$ MHz. The average length of segments used to discretize the whole computational domain is 0.033 m, which corresponds to $\lambda/20$, where $\lambda$ is the wavelength inside dielectric objects. 
	
	Fig. \ref{RCS} shows the RCS of the two dielectric objects obtained from the proposed hybrid SIE-PDE formulation with the reference RCS obtained from the FEM. It is easy to find that the results obtained from the two formulations for both objects agree well with each other. Therefore, the proposed hybrid SIE-PDE formulation can obtain accurate far-fields induced by penetrable objects with smooth and non-smooth boundaries.
	
	Then, we check the convergence property of the proposed SIE-PDE formulation in terms of mesh sizes. The relative error (RE) is defined as
	\begin{equation} \label{RE}
		\text{RE} = \frac{{\sum\limits_i {{{\left\| {\text{RCS}^\text{cal} - \text{RCS}^\text{ref}} \right\|}^2}} }}{{\sum\limits_i {{{\left\| \text{RCS}^\text{ref} \right\|}^2}} }},
	\end{equation}
	where $\text{RCS}^\text{cal}$ denotes results calculated from the proposed SIE-PDE formulation, and $\text{RCS}^\text{ref}$ is the reference results obtained from the FEM with fine enough meshes. As shown in Fig. \ref{convergence}, it can be found that as the mesh size decreases, the relative error of the proposed SIE-PDE formulation for both objects decreases. For the cylinder object, the RE can reach the order -4 when the average mesh size is 0.01. For the cuboid object, the relative error is 0.001 when the average mesh size is 0.01, which is slightly larger than that of the cylinder object. This is expected since the cuboid object has a non-smooth boundary and large fields variations exist near corners of the cuboid. However, they all show excellent agreement for these two formulations. Therefore, the proposed SIE-PDE formulation shows great convergence properties in terms of the mesh size.
\begin{figure*}[t]
		\begin{minipage}[h]{0.3\linewidth}
			\centerline{\includegraphics[scale=0.26]{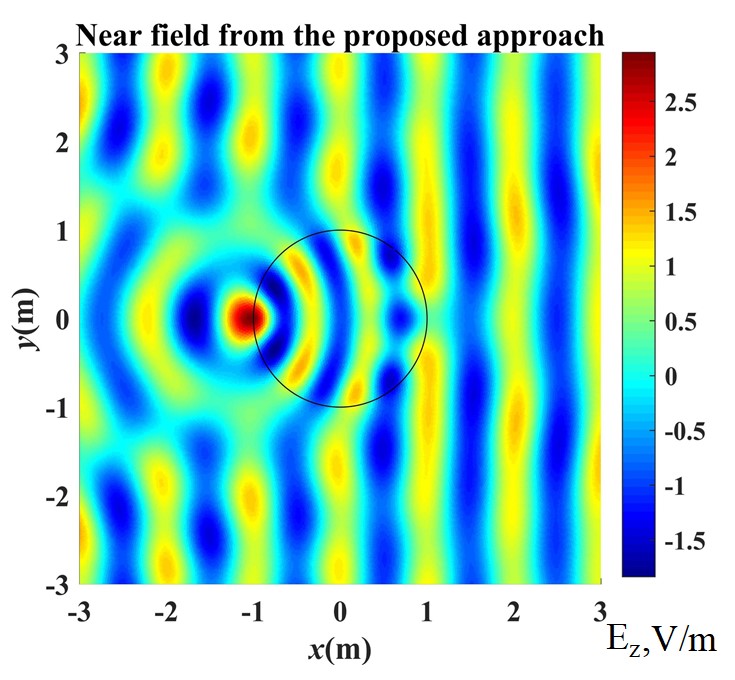}}
			\centerline{(a)}
		\end{minipage}
		\hfill
		\begin{minipage}[h]{0.3\linewidth}
			\centerline{\includegraphics[scale=0.26]{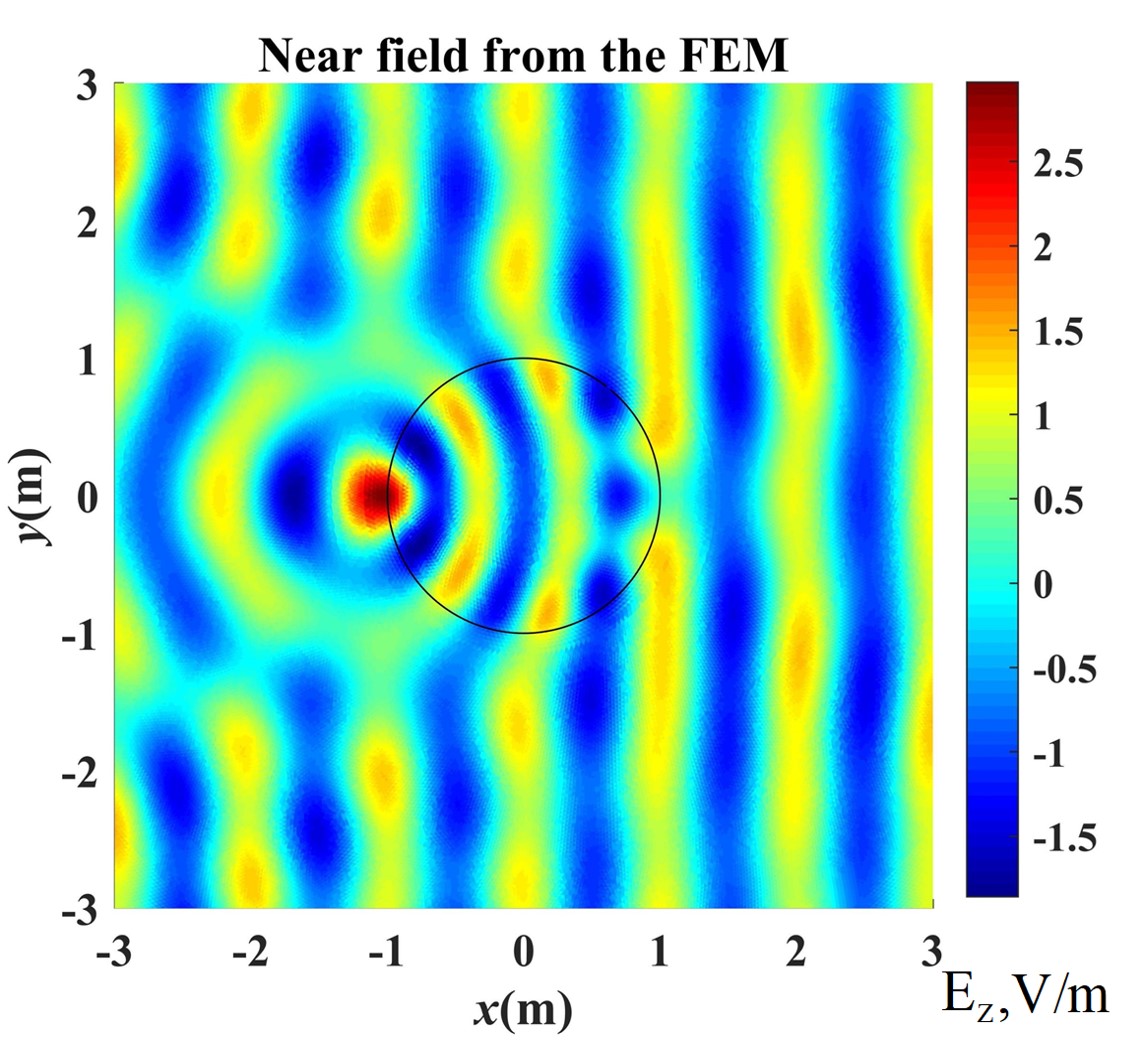}}
			\centerline{(b)}
		\end{minipage}
		\hfill
		\begin{minipage}[h]{0.3\linewidth}
			\centerline{\includegraphics[scale=0.26]{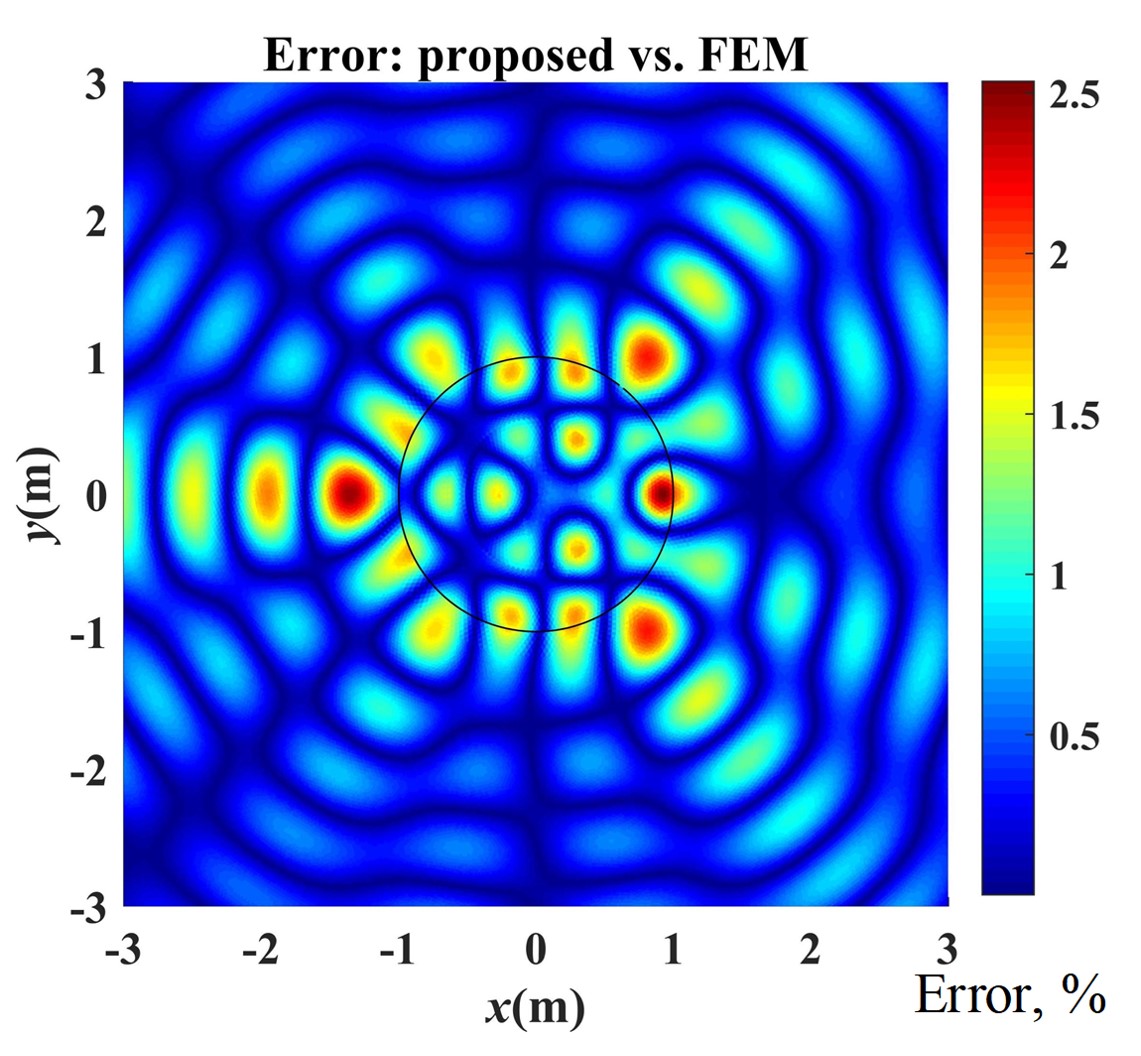}}
			\centerline{(c)}
		\end{minipage}
		\caption{(a) Near fields obtained from the proposed SIE-PDE formulation, (b) near fields obtained from the FEM, and (c) the relative error of near fields obtained from the proposed SIE-PDE formulation and the FEM.}
		\label{nearfield_circle}
	\end{figure*}
\begin{figure*}[t]
	\begin{minipage}[h]{0.3\linewidth}
		\centerline{\includegraphics[scale=0.26]{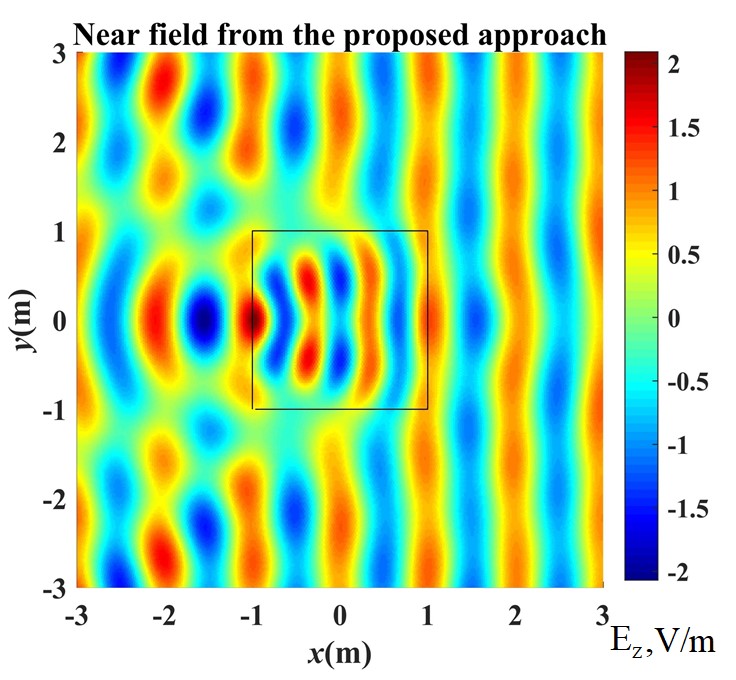}}
		\centerline{(a)}
	\end{minipage}
	\hfill
	\begin{minipage}[h]{0.3\linewidth}
		\centerline{\includegraphics[scale=0.26]{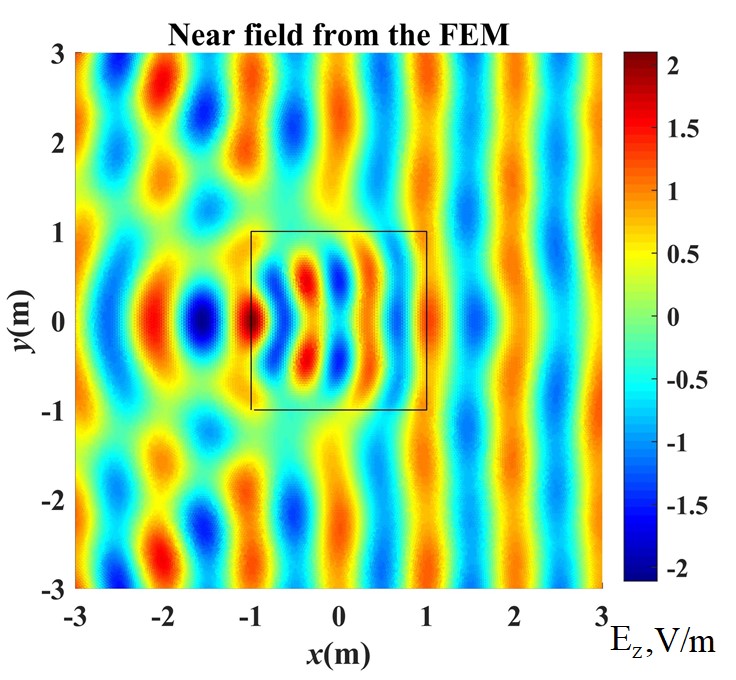}}
		\centerline{(b)}
	\end{minipage}
	\hfill
	\begin{minipage}[h]{0.3\linewidth}
		\centerline{\includegraphics[scale=0.26]{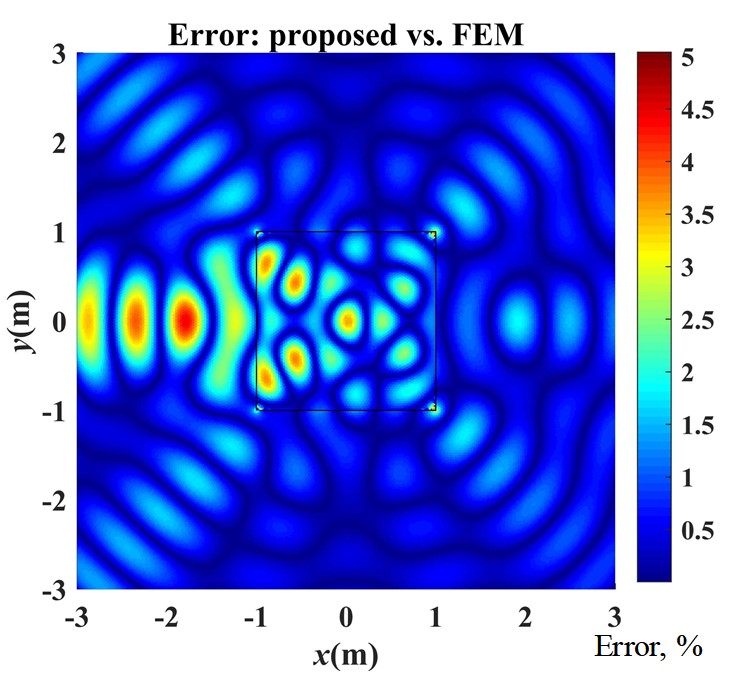}}
		\centerline{(c)}
	\end{minipage}
	\caption{(a) Near fields obtained from the proposed SIE-PDE formulation, (b) near fields obtained from the FEM, and (c) the relative error of near fields obtained from the proposed SIE-PDE formulation and the FEM.}
	\label{nearfield_rectangle}
\end{figure*}

	For further verification, we calculated the electric fields in the near region of both dielectric objects through the proposed SIE-PDE formulation and the FEM in Fig. \ref{nearfield_circle}(a) and (b), Fig. \ref{nearfield_rectangle}(a) and (b), respectively. It is easy to find that the field patterns obtained from the two formulations have no visible differences. To quantitatively measure the error between the two 
	formulations, the relative error defined as $\left| {{\text{E}^\text{cal}} - {\text{E}^\text{ref}}} \right|/\max \left| {{\text{E}^\text{ref}}} \right|$ is calculated, where ${\text{E}^\text{ref}}$, ${\text{E}^\text{cal}}$ denote the electric fields obtained from the FEM and the proposed SIE-PDE formulation, respectively, and $\max \left| {{\text{E}^\text{ref}}} \right|$ denotes the maximum magnitude of the reference electric field, which guarantees the relative error is well defined in the whole computational domain. 

	For the dielectric cylinder in Fig. \ref{nearfield_circle}(c), the relative error is less than 2$\%$ in most of the computational domain, and a slightly large relative error occurs near the boundary of the dielectric cylinder. It is expected since in those regions fields change sharply due to media discontinuity. For the dielectric cuboid in Fig. \ref{nearfield_rectangle}(c), the relative error in most regions is less than 3$\%$ and slightly large near the boundary of the cuboid object, especially at its four corners. However, the maximum of relative error is still around 4$\%$. Therefore, the proposed SIE-PDE formulation shows a quite good accuracy compared with the FEM. 
	\begin{figure}[t]
	\centering
	\includegraphics[width=0.37\textwidth]{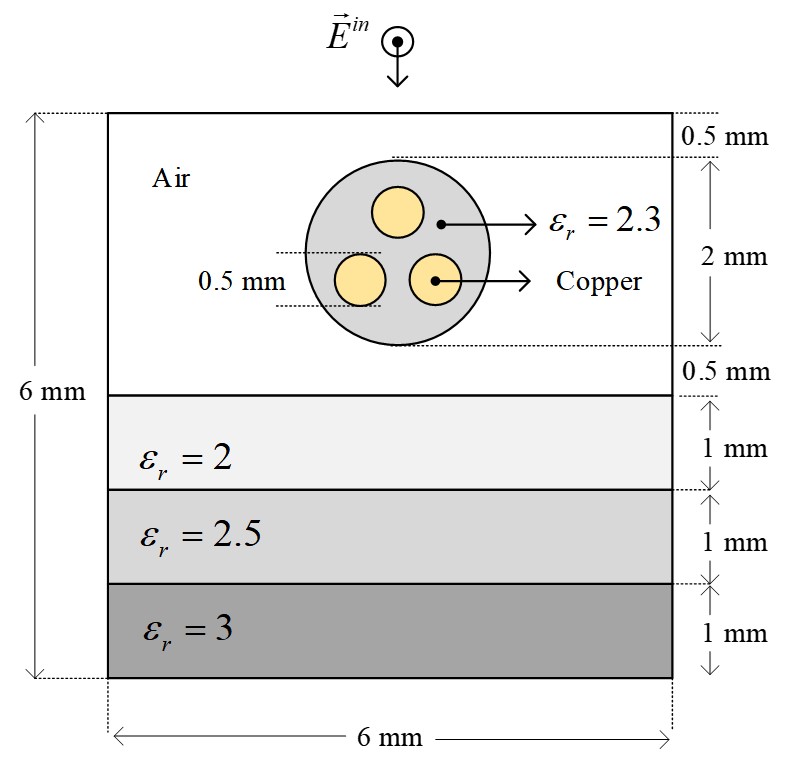}
	\caption{Geometrical configurations of the cable over three planar layered media.}
	\label{complex_geo} 
	\end{figure}
	
	To sum up, through the numerical results above, we find that the proposed SIE-PDE formulation can accurately calculate both near and far fields induced by dielectric structures with smooth and non-smooth boundaries, and shows great agreement with those obtained from the FEM. To further verify the computational efficiency and accuracy of the proposed SIE-PDE formulation, a more complex structure is considered in next subsection, and the computational costs will be discussed in detail. 

	\subsection{Computational Efficiency and Accuracy Verification of the Proposed Hybrid SIE-PDE Formulation}
	
\begin{figure}[t]
		\centering
		\includegraphics[width=0.44\textwidth]{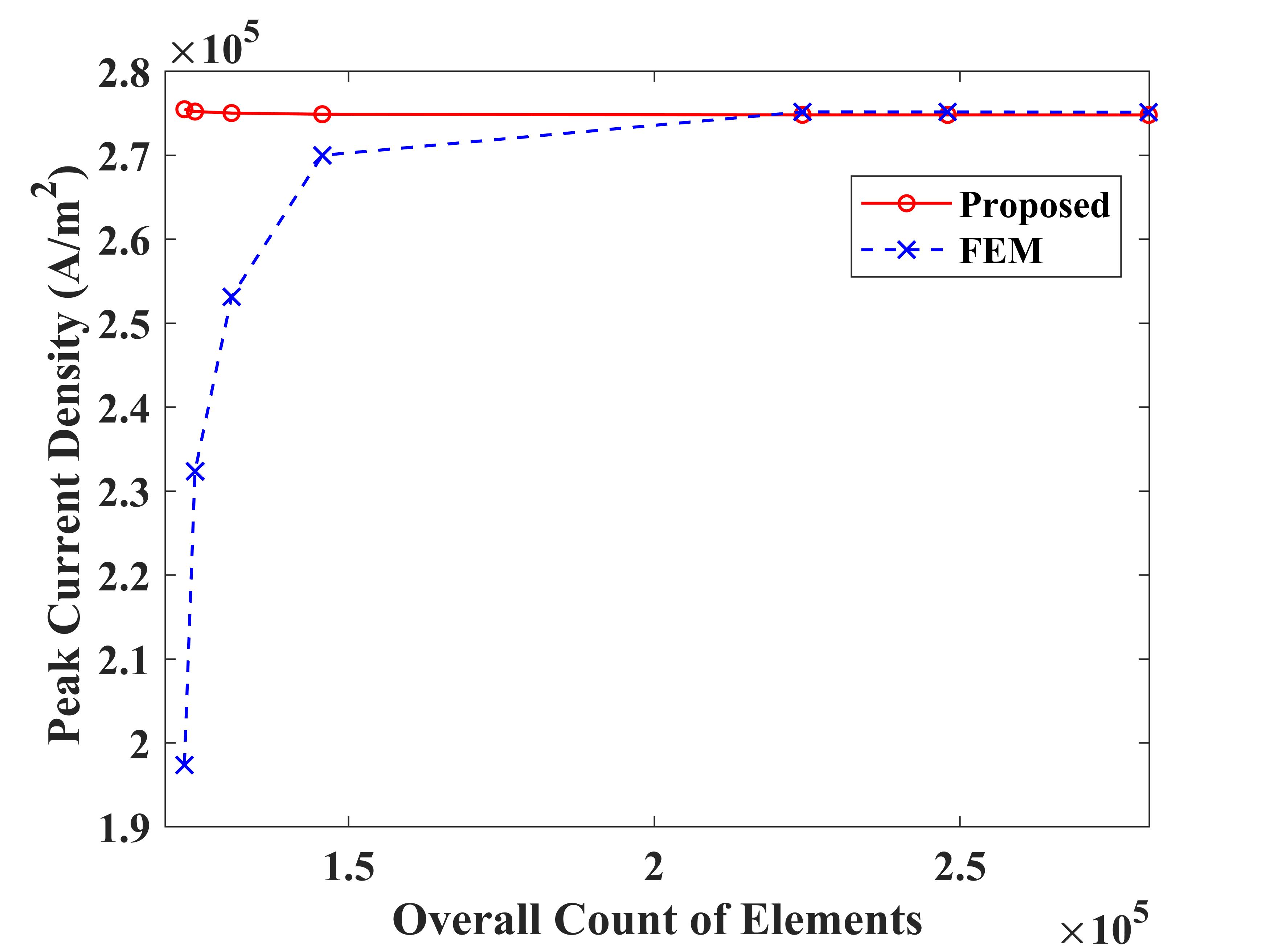}
		\caption{The peak current density obtained from the proposed SIE-PDE formulation and the FEM in terms of the count of elements.}
		\label{current_mesh} 
	\end{figure}

	Due to the skin effect, the current density crowds round the surface region of highly conductive media when the frequency is high. For example, the skin depth of copper at $f = 30$ MHz is approximately 11.9 $\mu m$. In the traditional volumetric mesh formulations, like the FEM, the volume integral equation (VIE) formulation, extremely fine meshes should be used to accurately calculate the current density. For example, to accurately capture the skin effect, almost 40,000 volumetric elements are used for only a $10\mu m \times 10\mu m \times 10\mu m$ copper interconnect at $f = 20$ GHz [\citen{AlQedra2010}]. Therefore, it is quite challenging to solve this problem in terms of CPU time and memory consumption. Compared with those volumetric meshes-based solvers, one obvious merit of the SIE formulations is that they can effectively address computational issues related to conductors and obtain accurate results since unknowns only reside on the boundary of highly conductive media. The proposed hybrid SIE-PDE formulation inherits this merit, and it also increases the flexibility to handle the inhomogeneous media, in which the SIE formulation cannot be used. To further verify the computational efficiency and accuracy of the proposed hybrid SIE-PDE formulation, a numerical example related to conductive media and skin effect calculation is considered. Its geometrical configurations are given in Fig. \ref{complex_geo}. One cable including three cylindrical copper conductors with the radius of 0.5 mm is placed above the planar layered media. The three cylindrical conductors are embedded in a dielectric sheath with the radius of 1 mm. The relative permittivity of dielectric cylinder is 2.3. Then, the cable system is placed 0.5 mm over the planar layered media with the relative permittivity of 1, 2, 2.5 and 3, respectively. A TM-polarized plane wave with $f = 30$ MHz incidents along the $y$-axis.
	\begin{figure*}[t]
	\begin{minipage}[h]{0.3\linewidth}
		\centerline{\includegraphics[scale=0.26]{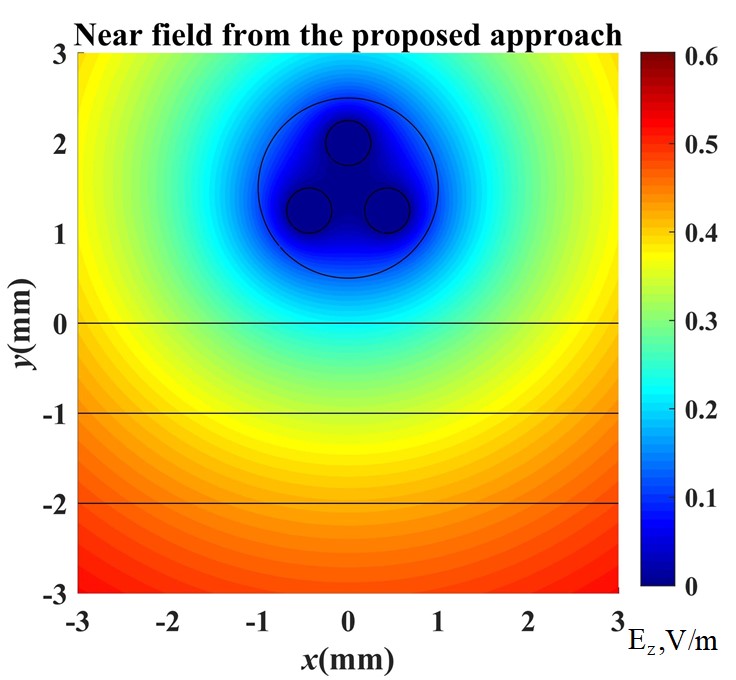}}
		\centerline{(a)}
	\end{minipage}
	\hfill
	\begin{minipage}[h]{0.3\linewidth}
		\centerline{\includegraphics[scale=0.26]{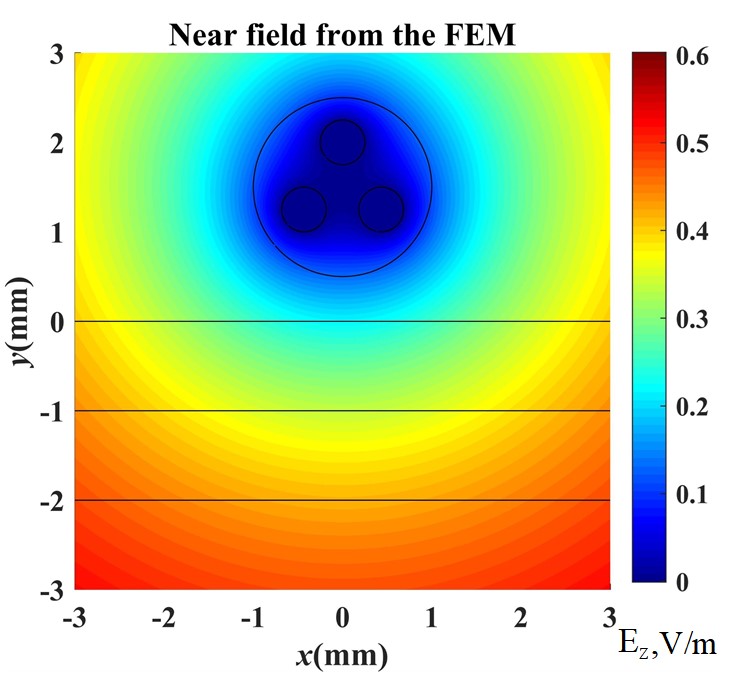}}
		\centerline{(b)}
	\end{minipage}
	\hfill
	\begin{minipage}[h]{0.3\linewidth}
		\centerline{\includegraphics[scale=0.26]{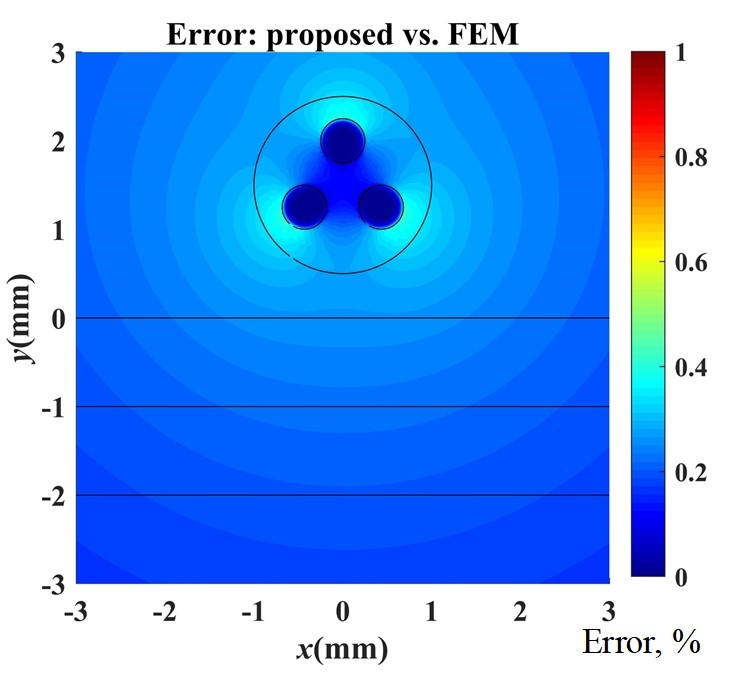}}
		\centerline{(c)}
	\end{minipage}
	\caption{(a) Near fields obtained from the proposed SIE-PDE formulation, (b) near fields obtained from the FEM, and (c) the relative error of near fields obtained from the proposed SIE-PDE formulation and the FEM.}
	\label{nearfield_mm}
\end{figure*}
	We calculated electric fields near the object through the proposed SIE-PDE formulation and the FEM, as shown in Fig. \ref{nearfield_mm}(a) and (b). It can be found that field patterns from the two formulations are almost the same. To quantitatively measure the error, we also calculated the relative error of near fields. As shown in Fig. \ref{nearfield_mm}(c), the relative error is less than 0.5$\%$ in most regions of the computational domain and slightly large on the boundary of the inner cylindrical conductors. In the current mesh configuration, both formulations use the same mesh size, 0.05 mm. However, by further calculating the current density in the conductors, the peak values of the current density obtained from the two formulations show significant differences. The values obtained from the proposed SIE-PDE formulation and the FEM are 275,492.7 $A/m^2$ and 197,378.4 $A/m^2$, respectively. Then, we checked the convergence property of the proposed SIE-PDE formulation and the FEM in terms of the current density in Fig. \ref{current_mesh}. The $x$-axis represents the overall count of elements in the whole computational domain, and the $y$-axis represents the peak amplitude of the current density. It can be found that as the count of elements increases, the current density obtained from the FEM gradually increases, and finally is convergent at around 275,000 $A/m^2$. However, from the proposed SIE-PDE formulation, the results get the convergent results quite fast with much fewer elements, which implies that the proposed SIE-PDE formulation can obtain the accurate current density with much coarser meshes than those with the FEM. When the count of elements reaches approximate $ 2.3\times {10^5}$, the FEM gets convergent results. For example, when the count of elements is 224,244, the peak values of the current density obtained from the proposed SIE-PDE formulation and the FEM are 274,827.9 $A/m^2$ and 275,180.0 $A/m^2$, respectively. The two methods have almost the same level of accuracy.	 
	
	\begin{figure}
	\begin{minipage}[t]{0.48\linewidth}
		\centerline{\includegraphics[scale=0.3]{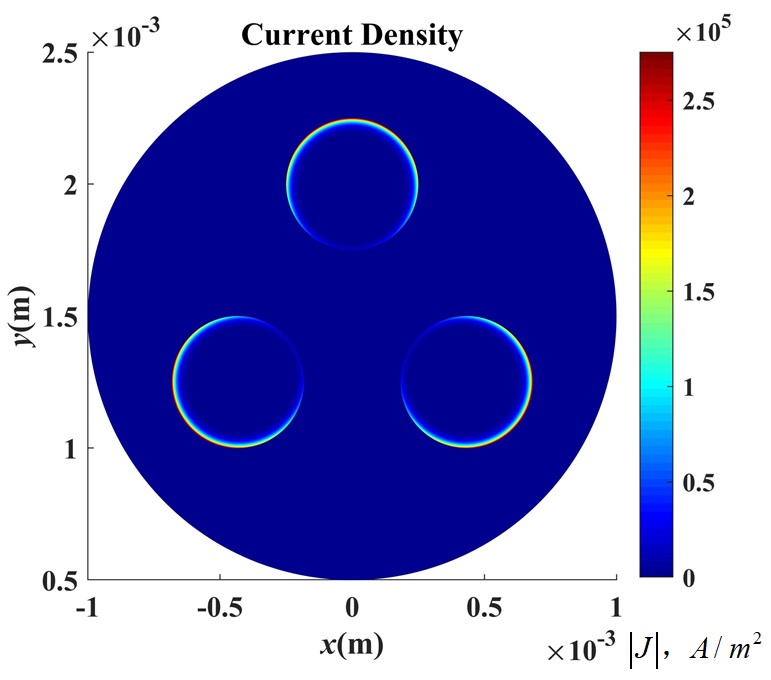}}
		\centerline{(a)}
	\end{minipage}
	\hfill
	\begin{minipage}[t]{0.48\linewidth}
		\centerline{\includegraphics[scale=0.3]{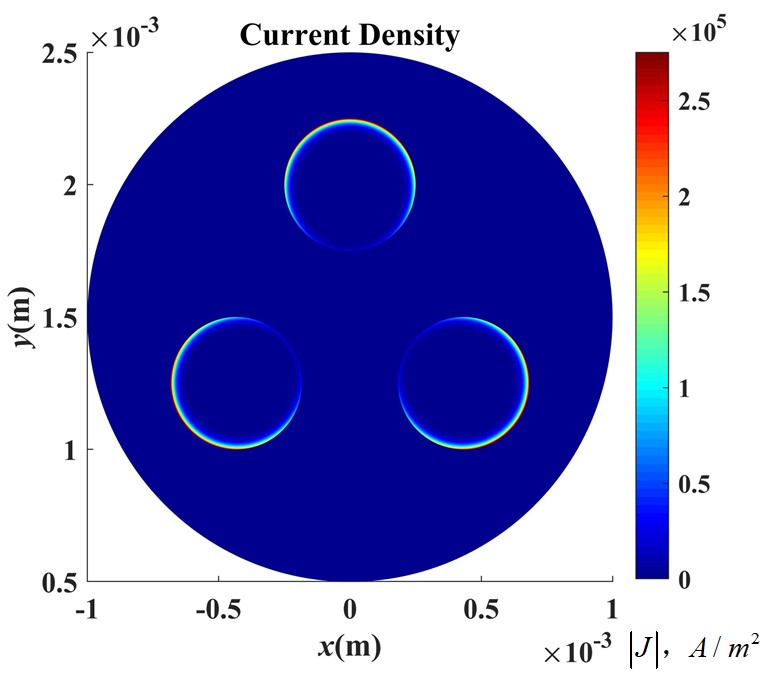}}
		\centerline{(b)}
	\end{minipage}
	\caption{The current density obtained from (a) the proposed SIE-PDE formulation, and (b) the FEM.}
	\label{current_density}
\end{figure}
	Fig. \ref{current_density}(a) and (b) show the current density obtained from the proposed SIE-PDE formulation and the FEM. It is easy to find that both formulations can accurately capture the skin effect in terms of the current density, and Table I shows the comparison of computational consumption of the two formulations when the current density with the same level of accuracy is obtained. The ratio in the last column of Table I is defined as the ratio of the values obtained from the proposed SIE-PDE formulation to those from the FEM. As shown in the second row of Table I, the peak value of the current density is 274,899.7 $A/m^2$ for the proposed SIE-PDE formulation compared with 274,965.9 $A/m^2$ for the FEM. Its ratio is 0.9997 and to this extent, the two formulations achieve almost the same level of accuracy. Then, we calculated the overall unknowns used for the two formulations. As shown in the third row of Table I, the proposed SIE-PDE formulation only needs 73,234 unknowns to solve this problem compared with 253,673 unknowns for the FEM. Its ratio is 0.29, which means that only quite a small fraction of amount of unknowns is required in the proposed SIE-PDE formulation. As is illustrated above, the SIE formulations can capture the skin effect inside the conductors fast and accurately without requirements for fine meshes. In this case, for the FEM the average length of segments is used as 0.005 mm for the inner cylindrical copper conductors and dielectric object and 0.05 mm for outer layered medium. However, in the proposed SIE-PDE formulation, the average length of segments is only 0.02 mm for the inside objects. Therefore, the overall count of unknowns is significantly reduced for the proposed formulation compared with the FEM. Moreover, as shown in the fourth row of Table I, memory consumption of the proposed SIE-PDE formulation is only 36$\%$ of that of the FEM. The time cost for the two formulations is shown from the fifth to eighth rows in Table I. The time cost includes generation of $\mathbb{Y}_s$, matrices filling and matrices solving. The time cost for the FEM requires 4,771.9 seconds to solve this problem. However, only 330.4 seconds are required for the proposed formulation, which only requires 7$\%$ of CPU time of the FEM. As shown in Table I, the time cost of matrices filling and solving linear equations is reduced compared with that of the FEM due to much fewer unknowns in the proposed SIE-PDE formulation. As the count of SIE domains increases, the performance improvement in terms of time cost and count of unknowns would be more significant.

\begin{table}[h]
	\centering
	\caption{Comparison of computational costs for the proposed SIE-PDE formulation and the FEM}\label{table1}
	\resizebox{8.8cm}{!}{
		\begin{threeparttable}[b]
			\begin{tabular}{l|c|c|c}
				\hline
				\hline
				\textbf{ } &\textbf{FEM}   & \textbf{Proposed} & {\textbf{Ratio}*}\cr 
				\hline
				\hline
				Peak Value of Current Density [$A/m^2$]        &{274,965.9}  & {274,899.7}         &{0.99}      \cr
				\hline
				Overall Counts of Unknowns      &{253,673}   & {73,234} & {0.36}    \\
				\hline
				\textbf{Total Time } [s] &4,771.9 &330.4   & 0.07  \\
				\hline
				Time for $\mathbb{Y}_s$ Generation [s]&{-}    &{232.0}  & {-}      \\
				\hline
				Time for Matrices Filling [s]      &{4,769.0} &          {97.7}  &           {0.02}       \\  
				\hline
				Time for Matrices Solving [s]         &{2.9}            &{0.7}     &  {0.24}   \\
				\hline 
				\hline 
			\end{tabular}
			\begin{tablenotes}
				\footnotesize
				\item[*]Ratio is defined as the ratio of the values obtained from the proposed SIE-PDE formulation to those from the FEM.
			\end{tablenotes}
		\end{threeparttable}
	}
\end{table}

	\subsection{Discussion}
	From the above numerical examples, the computational accuracy and efficiency of the proposed SIE-PDE formulation for 2D TM electromagnetic problems are detailed presented. It should be noted that the SIE formulation is used to solve the interior problems instead of truncating the whole computational domain. Therefore, we compared the proposed SIE-PDE formulation with the traditional FEM under the same absorbing boundary condition rather than other hybrid FEM/MoM formulations. 
	
	The proposed hybrid formulation can be extended into 3D scenarios and applied to solve the vector Helmholtz equation by exchanging the 2D nodal-based linear basis function in (\ref{E0_exp}) with 3D edge-based vector basis function [\citen{Jin2015FEM}, Ch, 8, pp. 420-435]. Then, to couple the SIE and PDE formulations on the interface, the Rao-Wilton-Glisson (RWG) basis function [\citen{Gibson2015}, Ch, 8, 258-260] and the 2D edge-based vector basis function will be used due to their relationship on the surface, ${\bf{n}} \times {{\bf{N}}_i} = {{\bf{f}}_n}$, where $\mathbf{N}_i$ denotes the vector basis function defined on the $i$th edge, and $\mathbf{f}_n$ denotes the RWG basis function defined on the same edge. In addition, the proposed SIE-PDE formulation can also support non-conformal meshes with appropriate interpolation methods.

	\section{Conclusion}
	In this paper, a novel and efficient hybrid SIE-PDE formulation for TM electromagnetic analysis without boundary condition requirement is developed. In this formulation, the whole computational domain is decomposed into two $\emph{overlapping}$ domains: the SIE and PDE domains. In the SIE domain, complex structures or computationally challenging media are replaced by the background medium. An equivalent model with only the surface equivalent electric current density on the enclosed boundary is derived to represent electromagnetic effects. The remaining computational domain and the background medium replaced domain consist of the PDE domain which can be formulated through the PDE formulation. Through carefully constructing the basis functions, the expansion of the SIE and PDE formulations are compatible on the boundary. Then, a hybrid SIE-PDE formulation is obtained. The proposed formulation shows many advantages over other existing hybrid techniques, such as no extra boundary condition requirements, flexibility and high efficiency in handling complex and computationally challenging structures. Numerical results demonstrate that the proposed SIE-PDE formulation can significantly improve the performances in terms of the numbers of unknowns and CPU time over the FEM, which shows great potential for electromagnetic analysis of complex, non-uniform media fulfilled objects. Additionally, it can be extended to solve 3D electromagnetic problems and support non-conformal meshes. Currently, the development of the proposed SIE-PDE formulation into 3D and non-conformal scenarios is in progress. We will report more results on this topic in the future.
		
%
%
%
%
%

	\section*{Appendix}
		In this paper, the Green’s function is expressed as $G\left( {\mathbf{r},\mathbf{r}'} \right) = -jH_0^{\left( 2 \right)}\left( {k\rho } \right)/4$, where $j = \sqrt {-1}$, $\rho  = \left| {\mathbf{r} - \mathbf{r}'} \right|$, $k$ denotes the wave number in the penetrable object, and $H_0^{\left( 2 \right)}\left(  \cdot  \right)$ denotes the zeroth-order Hankel function of the second kind. When $\left| {k\rho } \right| < 0.1$, the Green’s function and its gradient suffer from nearly singular and singular integrals. The first-order small variable approximation of $H_0^{\left( 2 \right)}\left(  \cdot  \right)$ and $H_1^{\left( 2 \right)}\left(  \cdot  \right)$ [\citen{Wiltonintegral1984}] is considered, which is expressed as 
		\begin{equation} \label{Hankel_0}
			H_0^{(2)}(k \rho)  \approx  \left[ {1 - j\frac{2}{\pi }\ln \left( {\frac{{\gamma k}}{2}} \right)} \right] - j\frac{2}{\pi }\ln \rho,
		\end{equation}
		\begin{equation} \label{Hankel_1}
			H_1^{(2)}(k\rho) \approx \frac{{k\rho}}{2} + \frac{{2j}}{{\pi k\rho}},
		\end{equation}
		where $\gamma  = 1.781$, $\rho  = \left| {\mathbf{r} - \mathbf{r}'} \right|$.
		
\begin{figure}[t]
			\centering
			\includegraphics[width=0.4\textwidth]{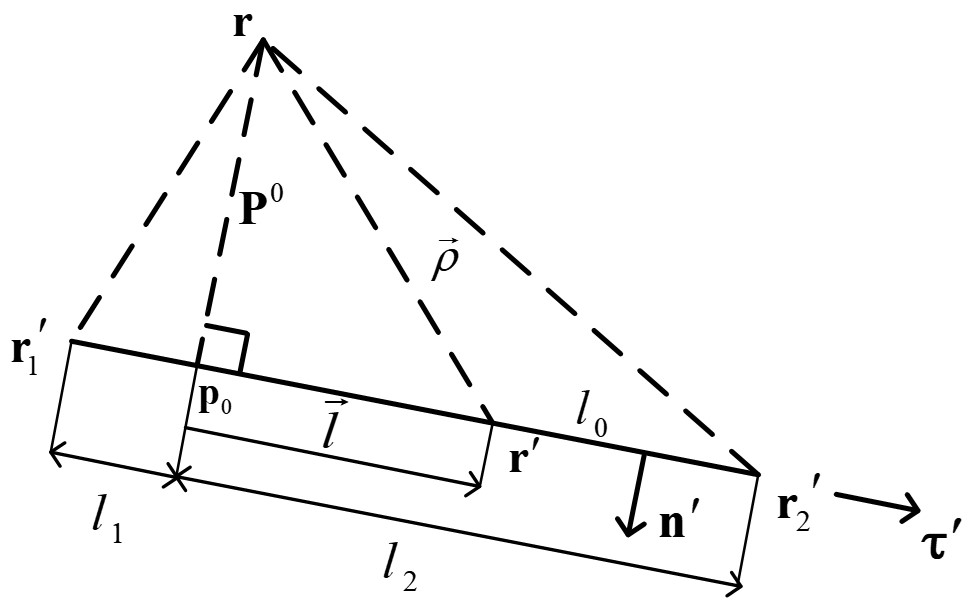}
			\caption{A general configuration to demonstrate the geometrical relationship between
				the observation point and the source segment. $\vec \rho  = \mathbf{r} - \mathbf{r}'$, $\mathbf{\tau}' = \frac{{{\mathbf{r}_2}^\prime  - {\mathbf{r}_1}^\prime }}{{\left| {{\mathbf{r}_2}^\prime  - {\mathbf{r}_1}^\prime } \right|}}$, $\vec l = \mathbf{\tau}' \cdot \left[ {\left( {\mathbf{r}' - \mathbf{r}} \right) \cdot \tau '} \right]$, ${l_0} = \left| {{\mathbf{r}_2}^\prime  - {\mathbf{r}_1}^\prime } \right|$, ${l_1} = \left( {{\mathbf{r}_1}^\prime  - \mathbf{r}} \right) \cdot \tau '$, ${l_2} = \left( {{\mathbf{r}_2}^\prime  - \mathbf{r}} \right) \cdot \tau '$, ${\mathbf{P}^0} = \left| {{\mathbf{p}_0} - \mathbf{r}} \right|$.}
			\label{green_all} 
		\end{figure}

		A general configuration to demonstrate the geometrical relationship between
		the observation point and the source segment is shown in Fig. \ref{green_all}. $\mathbf{r}$, $\mathbf{r}'$, ${{\mathbf{r}_1}^\prime }$, ${{\mathbf{r}_2}^\prime}$ denote the observation point, the source point, and two end points of the source segment, respectively, and ${{\mathbf{p}_0}}$ denotes the projection point of $\mathbf{r}$ on the source segment. As shown in (\ref{roogtop_fun}), there are two half rooftop basis functions defined on each line segment, $\frac{{\left| {\mathbf{r}' - {\mathbf{r}_1}^\prime } \right|}}{{{l_0}}}$ and $\frac{{\left| {{\mathbf{r}_2}^\prime  - \mathbf{r}'} \right|}}{{{l_0}}}$, respectively.
		Therefore, the analytical expression for (\ref{P1}) can be expressed as	
		\begin{align}
			\label{Hankel0_1_1}&{f_m}\left( {\mathbf{r}} \right) \cdot \int_{{\mathbf{r}_1}^\prime }^{{\mathbf{r}_2}^\prime } {j\omega \mu {G_1}\left( {\mathbf{r},\mathbf{r}'} \right) \cdot \frac{{\left| {\mathbf{r}' - {\mathbf{r}_1}^\prime } \right|}}{{{l_0}}}d\mathbf{r}'}  \\	
			&= {f_m}\left( {\mathbf{r}} \right) \cdot\left\{\frac{{\omega \mu }}{{4{l_0}}}\left[ {1 - j\frac{2}{\pi }\ln \left( {\frac{{\gamma k}}{2}} \right)} \right] \cdot {I_1} - j\frac{{\omega \mu }}{{2\pi {l_0}}} \cdot {I_2} \right. \left.+ j\frac{{\omega \mu }}{{2\pi {l_0}}}{l_1} \cdot {I_3} - \frac{{\omega \mu {l_1}}}{{4{l_0}}}\left[ {1 - j\frac{2}{\pi }\ln \left( {\frac{{\gamma k}}{2}} \right)} \right] \cdot {I_4} \right\}, \notag
		\end{align}
		\begin{align}
			\label{Hankel0_1_2}&{f_m}\left( {\mathbf{r}} \right) \cdot \int_{{\mathbf{r}_1}^\prime }^{{\mathbf{r}_2}^\prime } {j\omega \mu {G_1}\left( {\mathbf{r},\mathbf{r}'} \right) \cdot \frac{{\left| {{\mathbf{r}_2}^\prime - \mathbf{r}'} \right|}}{{{l_0}}}d\mathbf{r}'}  \\
			&= {f_m}\left( {\mathbf{r}} \right) \cdot \left\{  - \frac{{\omega \mu }}{{4{l_0}}}\left[ {1 - j\frac{2}{\pi }\ln \left( {\frac{{\gamma k}}{2}} \right)} \right] \cdot {I_1} + j\frac{{\omega \mu }}{{2\pi {l_0}}} \cdot {I_2} \right.\left.- j\frac{{\omega \mu }}{{2\pi {l_0}}}{l_2} \cdot {I_3} + \frac{{\omega \mu {l_2}}}{{4{l_0}}}\left[ {1 - j\frac{2}{\pi }\ln \left( {\frac{{\gamma k}}{2}} \right)} \right] \cdot {I_4} \right\}, \notag
		\end{align}
and the analytical expression for (\ref{U1}) can also be obtained as
		\begin{align}
			&{f_m}\left( {\mathbf{r}} \right) \cdot \int_{{\mathbf{r}_1}^\prime }^{{\mathbf{r}_2}^\prime } {k\frac{{\vec \rho  \cdot \mathbf{n}'}}{\rho }{G_1}^\prime \left( {\mathbf{r},\mathbf{r}'} \right)\frac{{\left| {\mathbf{r}' - {\mathbf{r}_1}^\prime } \right|}}{{{l_0}}}d\mathbf{r}'} \notag \\
			\label{Hankel1_1_1}
			&={f_m}\left( {\mathbf{r}} \right) \cdot \left\{ \frac{{j{k^2}}}{{8{l_0}}}{\mathbf{P}^0} \cdot {I_1} - \frac{{j{k^2}}}{{8{l_0}}}{\mathbf{P}^0}{l_1} \cdot {I_4}- \frac{{{\mathbf{P}^0}}}{{2\pi {l_0}}} \cdot {I_5}\right. \left. + \frac{{{\mathbf{P}^0}}}{{2\pi {l_0}}}{l_1} \cdot {I_6}\right\},
		\end{align}
		\begin{align}			
			&{f_m}\left( {\mathbf{r}} \right) \cdot \int_{{\mathbf{r}_1}^\prime }^{{\mathbf{r}_2}^\prime } {k\frac{{\vec \rho  \cdot \mathbf{n}'}}{\rho }{G_1}^\prime \left( {\mathbf{r},\mathbf{r}'} \right)\frac{{\left| {{\mathbf{r}_2}^\prime - \mathbf{r}'} \right|}}{{{l_0}}}d\mathbf{r}'}\notag \\ 
			\label{Hankel1_1_2}
			&= {f_m}\left( {\mathbf{r}} \right) \cdot \left\{ - \frac{{j{k^2}}}{{8{l_0}}}{\mathbf{P}^0} \cdot {I_1} + \frac{{j{k^2}}}{{8{l_0}}}{\mathbf{P}^0}{l_2} \cdot {I_4}+ \frac{{{\mathbf{P}^0}}}{{2\pi {l_0}}} \cdot {I_5} \right. \left. - \frac{{{\mathbf{P}^0}}}{{2\pi {l_0}}}{l_2} \cdot {I_6}\right\}.
		\end{align}

According to [\citen{Wiltonintegral1984}][\citen{StewartCal2007}], the six identities defined in (\ref{Hankel0_1_1})-(\ref{Hankel1_1_2}) are given as below.
		\begin{align} 
			\label{int_1}
			{I_1} =& \int_{{l_i}} {\left| {\vec l} \right|} d \mathbf{r}' = \frac{1}{2}\left[ {{{\left( {{l_2}} \right)}^2} - {{\left( {{l_1}} \right)}^2}} \right],\quad\\
			\label{int_2}
			{I_2} =& \int_{{l_i}} {\vec l}  \cdot \ln \left| {\mathbf{r} - \mathbf{r}'} \right|d\mathbf{r}' = \frac{{\vec \tau }}{2}\left\{ {{\left| {\mathbf{r} - {\mathbf{r}_2}^\prime } \right|}^2}\ln \left| {\mathbf{r} - {\mathbf{r}_2}^\prime} \right|- {{\left| {\mathbf{r} - {\mathbf{r}_1}^\prime } \right|}^2}\ln \left| {\mathbf{r} - {\mathbf{r}_1}^\prime } \right| - \frac{1}{2}\left[ {{{\left( {{l_2}} \right)}^2} - {{\left( {{l_1}} \right)}^2}} \right] \right\},\\
			\label{int_3}
			{I_3} =& \int_{{l_i}} {\ln \left| {\mathbf{r} - \mathbf{r}'} \right|d\mathbf{r}'}  = {l_2}\ln \left| {\mathbf{r} - {\mathbf{r}_2}^\prime } \right| - {l_1}\ln \left| {\mathbf{r} - {\mathbf{r}_1}^\prime } \right| + {\mathbf{P}^0}\left( {{{\tan }^{ - 1}}\frac{{{l_2}}}{{{\mathbf{P}^0}}} - {{\tan }^{ - 1}}\frac{{{l_1}}}{{{\mathbf{P}^0}}}} \right) - \left( {{l_2} - {l_1}} \right),\\
			\label{int_4}
			{I_4} =& \int_{{l_i}} 1 d\mathbf{r}' = {l_2} - {l_1},\\
			\label{int_5}
			{I_5} =& \int_{{l_i}}\! {\frac{{\left| {\vec l} \right|}}{{{{\left| {\mathbf{r} - \mathbf{r}'} \right|}^2}}}d\mathbf{r}'}\!  =\! \frac{1}{2}\left( {\ln {{\left| {\mathbf{r} - {\mathbf{r}_2}\!^\prime } \right|}^2}\! -\! \ln {{\left| {\mathbf{r} - {\mathbf{r}_1}\!^\prime }\! \right|}^2}} \right)\!,\\
			\label{int_6}
			{I_6} =& \int_{{l_i}}\!{\frac{1}{{{{\left| {\mathbf{r} - \mathbf{r}'} \right|}^2}}}d\mathbf{r}'} \! =\! \frac{1}{{{\mathbf{P}^0}}}\! \cdot\! \left[\!{{{\tan }^{ - 1}}\!\left( {\frac{{{l_2}}}{{{\mathbf{P}^0}}}} \right) \!-\! {{\tan }^{ - 1}}\!\left( {\frac{{{l_1}}}{{{\mathbf{P}^0}}}} \right)}\! \right]. \quad
			\end{align}


\begin{thebibliography}{00}	
		\bibitem{PengNon-DDM2013}	
		Z. Peng, K. Lim, and J. Lee, ``Nonconformal Domain Decomposition Methods for Solving Large Multiscale Electromagnetic Scattering Problems,'' {\it Proc. IEEE}, vol. 101, no. 2, pp. 298-319, Feb. 2013.	
		
		\bibitem{GedneyFEM/MoM1992}
		S. D. Gedney, J. Lee, and R. Mittra, ``A combined FEM/MoM approach to analyze the plane wave diffraction by arbitrary gratings,'' {\it IEEE Trans. Microw. Theory Tech.}, vol. 40, no. 2, pp. 363-370, Feb. 1992.
		
		\bibitem{TeiFDTD-FETD2007}
		F. L.Teixeira, ``FDTD/FETD methods: a review on some recent advances and selected applications,'' {\it J. Microwaves Optoelectr.}, vol. 6, no. 1, pp. 83-95, Aug. 2007.
		
		\bibitem{YeungFDTD-FETD1999}
		M. S. Yeung, ``Application of the hybrid FDTD–FETD method to dispersive materials,'' {\it Microw. Opt. Tech. Lett.}, vol. 23, no. 4, pp. 238-242, Oct. 1999.
		
		\bibitem{SunFDTD-SETD2019}
		Q. Sun, R. Zhang, Q. Zhan, and Q. Liu, ``3-D Implicit–Explicit Hybrid Finite Difference/Spectral Element/Finite Element Time Domain Method Without a Buffer Zone,'' {\it IEEE Trans. Antennas Propag.}, vol. 67, no. 8, pp. 5469-5476, Aug. 2019.
		
		\bibitem{FanUS-FETD2020}
		K. Fan, B. Wei, and X. He, ``The implement of the hybrid meshes in US-FETD by using Domain decomposition method,'' in {\it Proc. IEEE MTT-S Int. Conf. Numer. Electromagn. Multiphys. Modeling Optim. (NEMO)}, Dec. 2020, pp. 1-3.
		
		\bibitem{ZhangFETD-GSM2017}
		K. Zhang, C. Wang, and J. Jin, ``A hybrid FETD-GSM algorithm for broadband full-wave modeling of resonant waveguide devices,'' {\it IEEE Trans. Microw. Theory Tech.}, vol. 65, no. 9, pp. 3147-3158, Sept. 2017.
		
		\bibitem{AydoMoM/MM2019}
		A. Aydoğan, ``A hybrid MoM/MM method for fast analysis of E-plane dielectric loaded waveguides,'' {\it AEU-Int. J. Electron. C.}, vol. 100, pp. 9-15, Dec. 2018.
		
		\bibitem{CabaMoM/MM2011}
		E. D. Caballero, H. Esteban, A. Belenguer, V. E. Boria, J. V. Morro, and J. Cascon, ``Efficient design of substrate integrated waveguide filters using a hybrid MoM/MM analysis method and efficient rectangular waveguide design tools,'' in {\it IEEE Int. Electromagn. Adv. Appl. Conf.}, Turin, Italy, 2011, pp. 456–459.
		
		\bibitem{ChenMoM-PO2007}
		M. Chen, Y. Zhang, X. Zhao, and C. Liang, ``Analysis of antenna around NURBS surface with hybrid MoM-PO technique,'' {\it IEEE Trans. Antennas Propag.}, vol. 55, no. 2, pp. 407-413, Feb. 2007.
		
		\bibitem{LiMoM-PO2017}
		L. Xiao, X. Wang, B. Wang, G. Zheng, and P. Chen, ``An efficient hybrid method of iterative MoM-PO and equivalent dipole-moment for scattering from electrically large objects,'' {\it IEEE Antennas Wireless Propag. Lett.}, vol. 16, pp. 1723-1726, Feb. 2017.
		
		\bibitem{GongMoM-PO2006}
		Z. Gong, B. Xiao, G. Zhu, and H. Ke, ``Improvements to the hybrid MoM-PO technique for scattering of plane wave by an infinite wedge,'' {\it IEEE Trans. Antennas Propag.}, vol. 54, no. 1, pp. 251-255, Jan. 2006.
		
		\bibitem{IlicFEM-MoM_ante_2009}
		M. M. Ilic, M. Djordjevic, A. Z. Ilic, and B. M. Notaro, ``Higher order hybrid FEM-MoM technique for analysis of antennas and scatterers,'' {\it IEEE Trans. Antennas Propag.}, vol. 57, no. 5, pp. 1452-1460, May 2009.
		
		\bibitem{HoppeFEM/MOM1994}
		D. J. Hoppe, L. W. Epp, and J. Lee, ``A hybrid symmetric FEM/MOM formulation applied to scattering by inhomogeneous bodies of revolution,'' {\it IEEE Trans. Antennas Propag.}, vol. 42, no. 6, pp. 798-805, Jun. 1994.
		
		\bibitem{XuFEM/MoM2016}
		R. Xu, L. Guo, H. He, and W. Liu, ``A hybrid FEM/MoM technique for 3-D electromagnetic scattering from a dielectric object above a conductive rough surface,'' {\it IEEE Geosci. Remote Sens. Lett.}, vol. 13, no. 3, pp. 314-318, Mar. 2016. 
		
		\bibitem{RenFEM/MoM2016}
		Y. Ren, Q. Liu, and Y. Chen, ``A hybrid FEM/MoM method for 3-D electromagnetic scattering in layered medium,'' {\it IEEE Trans. Antennas Propag.}, vol. 64, no. 8, pp. 3487-3495, Aug. 2016.
		
		\bibitem{IlicFEM-MoM_3D_2009}
		M. M. Ilic and B. M. Notaros, ``Higher order FEM-MoM domain decomposition for 3-D electromagnetic analysis,'' {\it IEEE Antennas Wireless Propag. Lett.}, vol. 8, pp. 970-973, Aug. 2009.
		
		\bibitem{JiFEM/MoM2002}
		Y. Ji and T. H. Hubing,``On the modeling of a gapped power-bus structure using a hybrid FEM/MoM approach,'' {\it IEEE Trans. Electromagn. Compat.}, vol. 44, no. 4, pp. 566-569, Nov. 2002. 
		
		\bibitem{TapMoM-UTD2005}
		K. Tap, T. Lertwiriyaprapa, P. H. Pathak, and K. Sertel, ``A hybrid MoM-UTD analysis of the coupling between large multiple arrays on a large platform,'' in {\it Proc. IEEE Antennas Propag. Soc. Int. Symp.}, 2006, vol. 4A, pp. 175–178.
		
		\bibitem{LiuMoM-UTD2010}
		Z. Liu, J. Yang, and C. Liang, ``The hybrid higher order MoM‐UTD formulation for electromagnetic radiation problems,'' {\it Microwave Opt. Technol. Lett.}, vol. 52, no. 5, pp. 1087-1091, Mar. 2010.
		
		\bibitem{EibertFE/BI1999}
		T. F. Eibert, J. L. Volakis, D. R. Wilton, and D. R. Jackson, ``Hybrid FE/BI modeling of 3-D doubly periodic structures utilizing triangular prismatic elements and an MPIE formulation accelerated by the Ewald transformation,'' {\it IEEE Trans. Antennas Propag.}, vol. 47, no. 5, pp. 843-850, May 1999.
		
		\bibitem{ShengFE/BI2002}
		X. Q. Sheng and E. K. Yung, ``Implementation and experiments of a hybrid algorithm of the MLFMA-enhanced FE-BI method for open-region inhomogeneous electromagnetic problems,'' {\it IEEE Trans. Antennas Propag.}, vol. 50, no. 2, pp. 163-167, Feb. 2002.
		
		\bibitem{MeyerFE-BI1994}
		F. J. C. Meyer, ``The Two-Dimensional finite element/boundary element method in electromagnetics: formulation applications error estimates and mesh adaptive procedures,'' Ph.D. dissertation, Univ. Stellenbosch, Stellenbosch, South Africa, 1994.
		
		\bibitem{ZhaoFEM-BEM2006}
		K. Zhao, M. N. Vouvakis, and J. Lee, ``Solving electromagnetic problems using a novel symmetric FEM-BEM approach,'' {\it IEEE Trans. Magn.}, vol. 42, no. 4, pp. 583-586, Apr. 2006.
		
		\bibitem{YangFE-BI-MLFMA2013}
		M. Yang, H. Gao, and X. Sheng, ``Parallel domain-decomposition-based algorithm of hybrid FE-BI-MLFMA method for 3-D scattering by large inhomogeneous objects,'' {\it IEEE Trans. Antennas Propag.}, vol. 61, no. 9, pp. 4675-4684, Sep. 2013.
		
		\bibitem{VouvakisFEM-IE2004}
		M. N. Vouvakis, S. Lee, K. Zhao, and J. Lee, ``A symmetric FEM-IE formulation with a single-level IE-QR algorithm for solving electromagnetic radiation and scattering problems,'' {\it IEEE Trans. Antennas Propag.}, vol. 52, no. 11, pp. 3060-3070, Nov. 2004.
		
		\bibitem{Guanmultisolver2017}
		J. Guan, S. Yan, and J. Jin, ``A multi-solver scheme based on combined field integral equations for electromagnetic modeling of highly complex objects,'' {\it IEEE Trans. Antennas Propag.}, vol. 65, no. 3, pp. 1236-1247, Mar. 2017. 
		
		\bibitem{Guanmultisolver2016}
		J. Guan, S. Yan, and J. Jin, ``A multisolver scheme based on robin transmission conditions for electromagnetic modeling of highly complex objects,'' {\it IEEE Trans. Antennas Propag.}, vol. 64, no. 12, pp. 5345-5358, Dec. 2016.
		
		\bibitem{DodigFEM-BEM2021}
		H. Dodig, D. Poljak, and M. Cvetković, ``On the edge element boundary element method/finite element method coupling for time harmonic electromagnetic scattering problems,'' {\it Int. J. Numer. Meth. Engng.}, vol. 122, no. 14, pp. 3613–3652, Apr. 2021.
		
		\bibitem{Jin2015FEM}
		J. Jin, {\it The Finite Element Method in Electromagnetics}. Hoboken, NJ, USA: Wiley, 2014.
		
		\bibitem{LeeDDM2005}
		S. Lee, M. N. Vouvakis, and J. Lee, ``A non-overlapping domain decomposition method with non-matching grids for modeling large finite antenna arrays,'' {\it J. Comput. Phys.}, vol. 203, no. 1, pp. 1-21, Feb. 2005.
		
		\bibitem{StupfelDDM2000}
		B. Stupfel and Martine Mognot, ``A domain decomposition method for the vector wave equation,'' {\it IEEE Trans. Antennas Propag.}, vol. 48, no. 5, pp. 653-660, May 2000.
		
		\bibitem{PengIEDDM2011}
		Z. Peng, X. Wang, and J. Lee, ``Integral equation based domain decomposition method for solving electromagnetic wave scattering from non-penetrable objects,'' {\it IEEE Trans. Antennas Propag.}, vol. 59, no. 9, pp. 3328-3338, Sep. 2011. 
		
		\bibitem{ZhaoDDM2008}
		K. Zhao, V. Rawat, and J. Lee, ``A domain decomposition method for electromagnetic radiation and scattering analysis of multi-target problems,'' {\it IEEE Trans. Antennas Propag.}, vol. 56, no. 8, pp. 2211-2221, Aug. 2008.
		
		\bibitem{LueDDM2018}
		J. Lu, Y. Chen, D. Li, and J. Lee, ``An embedded domain decomposition method for electromagnetic modeling and design,'' {\it IEEE Trans. Antennas Propag.}, vol. 67, no. 1, pp. 309-323, Jan. 2018. 
		
		\bibitem{LuSI-DDM2018}
		J. Lu and J. Lee, ``Signal integrity analysis of integrated circuits by using embedded domain decomposition method,'' {\it IEEE Trans. Microw. Theory Tech.}, vol. 66, no. 12, pp. 5369-5382, Dec. 2018.
		
		\bibitem{SunSIE-PDE2021}
		A. Sun and S. Yang, ``A hybrid SIE-PDE formulation without additional boundary conditions for electromagnetic analysis,'' accepted by {\it Proc. IEEE Antennas Propag. Soc. Int. Symp.}, 2021.
		
		\bibitem{LOVE}
		C. Balanis, \textit{Antenna Theory: Analysis and Design}. 3rd ed. Wiley, 2005.
		
		\bibitem{Zhou2021embedded}
		X. Zhou, Z. Zhu, and S. Yang, ``Towards a unified approach to electromagnetic analysis of objects embedded in multilayers,'' {\it J. Comput. Phys.}, vol. 427, pp. 110073, Feb. 2021.
		
		\bibitem{Okoshi1985}
		T. Okoshi, \textit{Planar Circuits for Microwaves and Lightwaves}, Springer Berlin Heidelberg, 1985.
		
		\bibitem{Patel2017SS-SIE}
		U. R. Patel, P. Triverio, and S. V. Hum, ``A single-source surface integral equation formulation for composite dielectric objects,'' in {\it Proc. IEEE Int. Symp. Antennas Propag.}, Jul. 2017, pp. 1453–1454.
		
		\bibitem{Knockaert2008}
		L. Knockaert, D. De Zutter, G. Lippens, and H. Rogier, ``On the Schur complement form of the Dirichlet-to-Neumann operator,'' {\it Wave Motion}, vol. 45, no. 3, pp. 309–324, Jan. 2008. 
		
		\bibitem{Patel2016contour}
		U. R. Patel and P. Triverio, ``Skin effect modeling in conductors of arbitrary shape through a surface admittance operator and the contour integral method,'' {\it IEEE Trans. Microw. Theory Tech.}, vol. 64, no. 9, pp. 2708-2717, Aug. 2016. 
		
		\bibitem{Zhou2021SS-SIE}
		X. Zhou, Z. Zhu, and S. Yang, ``Formulation of Single-Source Surface Integral Equation for Electromagnetic Analysis of 2D Partially Connected Penetrable Objects,'' {\it IEEE J. Multiscale and Multiphys. Comput. Techn.}, vol. 6, pp. 16-23, Jan. 2021.
		
		\bibitem{AlQedra2010}
		M. A. I. AlQedra, \textit{Surface impedance formulation for electric field integral equation in magneto-quasistatic and full-wave boundary element models of interconnects}, 2010.
		
		\bibitem{Wiltonintegral1984}
		D. Wilton, S. Rao, A. Glisson, D. Schaubert, O. Al-Bundak, and C. Butler, ``Potential integrals for uniform and linear source distribution on polygonal and polyhedral domains,'' {\it IEEE Trans. Antennas Propag.}, vol. 32, no. 3, pp. 276-281, Mar. 1984. 
		
		\bibitem{StewartCal2007}
		J. Stewart, \textit{Calculus, Brooks Cole}, 2007.
		
		\bibitem{Gibson2015}
		W. Gibson, \textit{The Method of Moments in Electromagnetics}, CRC press, 2015.
		
		%
		%
		%
		%
		%
		%
		
		
	\end{thebibliography}
\end{document}